# Toward optimal multistep forecasts in non-stationary autoregressions

CHING-KANG ING[1], JIN-LUNG LIN[2] and SHU-HUI YU[3]

[1]*Institute of Statistical Science, Academia Sinica, No. 128, Sec. 2, Academia Rd., Taipei 115, Taiwan, R.O.C. E-mail:* cking@stat.sinica.edu.tw
[2]*Department of Finance, National Dong Hwa University, No. 1, Sec. 2, Da Hsueh Rd., Shoufeng, Hualien 97401, Taiwan, R.O.C. E-mail:* Jlin@mail.ndhu.edu.tw
[3]*Institute of Statistics, National University of Kaohsiung, No. 700, Kaohsiung University Rd., Kaohsiung 81148, Taiwan, R.O.C. E-mail:* shuhui@nuk.edu.tw

This paper investigates multistep prediction errors for non-stationary autoregressive processes with both model order and true parameters unknown. We give asymptotic expressions for the multistep mean squared prediction errors and accumulated prediction errors of two important methods, plug-in and direct prediction. These expressions not only characterize how the prediction errors are influenced by the model orders, prediction methods, values of parameters and unit roots, but also inspire us to construct some new predictor selection criteria that can ultimately choose the best combination of the model order and prediction method with probability 1. Finally, simulation analysis confirms the satisfactory finite sample performance of the newly proposed criteria.

## 1. Introduction

Forecasting theory for stationary series with known true parameters is well studied but not much is known about the case for non-stationary models with estimated parameters. To fill the gap, this paper investigates multistep prediction errors for autoregressive (AR) processes with unit root. The plug-in and direct predictors are the two most frequently used multistep prediction methods and comparing their relative performance has become a major issue in forecast theory. In the case of squared error losses, the plug-in predictor is obtained from repeatedly using the fitted (by least squares) AR model with an unknown future value replaced by their own forecasts and the direct predictor is obtained by estimating the coefficient vector in the associated multistep prediction formula directly by linear least squares (see (1.2) and (1.3) below). Recently, many informative guidelines have been proposed to choose between these two methods in various time series models; see Findley [5, 6], Tiao and Tsay [20], Lin and Tsay [16], Ing [9, 10], Chevillon and Hendry [3] and Lin and Wei [17], among many others. However, a theoretical resolution to the problem of how to select the optimal multistep predictor in non-stationary time series still seems to be lacking, at least when the estimation uncertainty is taken into account. In this paper, we have developed and rigorously analyzed the theoretical properties of some predictor selection criteria to choose the model order and prediction method simultaneously.

Assume that observations $x_1, \ldots, x_n$ are generated from a unit root AR model,

$$x_{t+1} = \sum_{i=1}^{p+1} a_i x_{t+1-i} + \varepsilon_{t+1}, \tag{1.1}$$

where $0 \leq p < \infty$ is unknown, $a_{p+1} \neq 0$, $\varepsilon_t$'s are white noises with zero means and common variance $\sigma^2$ and the characteristic polynomial

$$A(z) = 1 - a_1 z - \cdots - a_p z^p - a_{p+1} z^{p+1}$$
$$= (1 - z)(1 - \alpha_1 z - \cdots - \alpha_p z^p),$$

with $\alpha(z) = (1 - \alpha_1 z - \cdots - \alpha_p z^p) \neq 0$ for all $|z| \leq 1$. $x_t$ is called stationary or stable if all roots of $A$ are outside the unit circle and unstable or non-stationary if some roots of $A$ are on the unit circle. For the sake of convenience, the initial conditions are set to $x_t = 0$ for all $t < 0$. To predict $x_{n+h}, h \geq 1$, based on $x_1, \ldots, x_n$ and a working model AR($k$), one may use the plug-in predictor, $\hat{x}_{n+h}(k)$, or direct predictor, $\check{x}_{n+h}(k)$, where

$$\hat{x}_{n+h}(k) = \mathbf{x}'_n(k)\hat{\mathbf{a}}_n(h, k), \tag{1.2}$$

and

$$\check{x}_{n+h}(k) = \mathbf{x}'_n(k)\check{\mathbf{a}}_n(h, k), \tag{1.3}$$

with $\mathbf{x}_j(k) = (x_j, \ldots, x_{j-k+1})'$ being the regressor vector and $\hat{\mathbf{a}}_n(h, k)$ and $\check{\mathbf{a}}_n(h, k)$ being plug-in and direct estimators, respectively. Note that

$$\left\{ \sum_{j=k}^{i-1} \mathbf{x}_j(k)\mathbf{x}'_j(k) \right\} \hat{\mathbf{a}}_i(1, k) = \sum_{j=k}^{i-1} \mathbf{x}_j(k) x_{j+1},$$

$$\left\{ \sum_{j=k}^{i-h} \mathbf{x}_j(k)\mathbf{x}'_j(k) \right\} \check{\mathbf{a}}_i(h, k) = \sum_{j=k}^{i-h} \mathbf{x}_j(k) x_{j+h},$$

and $\hat{\mathbf{a}}_i(h, k) = \hat{A}_i^{h-1}(k)\hat{\mathbf{a}}_i(1, k)$, with $\hat{A}_i^0(k) = I_k$,

$$\hat{A}_i(k) = \left( \hat{\mathbf{a}}_i(1, k) \; \middle| \; \frac{I_{k-1}}{\mathbf{0}'_{k-1}} \right),$$

and $I_m$ and $\mathbf{0}_m$, respectively, denoting an identity matrix and a vector of zeros of dimension $m$. To assess the prediction performance of $\hat{x}_{n+h}(k)$ and $\breve{x}_{n+h}(k)$, we consider their mean squared prediction errors (MSPEs),

$$\text{MSPE}P_{n,h}(k) = E(x_{n+h} - \hat{x}_{n+h}(k))^2$$

and

$$\text{MSPE}D_{n,h}(k) = E(x_{n+h} - \breve{x}_{n+h}(k))^2.$$

Theoretical investigations of $\text{MSPE}P_{n,h}(k)$ (or $\text{MSPE}D_{n,h}(k)$) in non-stationary AR models date back at least to Fuller and Hasza [7]. When $k \geq p+1$, an argument similar to that used in their Theorem 3.1 yields the following asymptotic expressions:

$$\text{MSPE}P_{n,h}(k) = \sigma_h^2 + E\{R_{P,n}(k)\} \tag{1.4}$$

and

$$\text{MSPE}D_{n,h}(k) = \sigma_h^2 + E\{R_{D,n}(k)\}, \tag{1.5}$$

where $R_{P,n}(k) = \text{O}_p(n^{-1})$, $R_{D,n}(k) = \text{O}_p(n^{-1})$ and $\sigma_h^2 = E(\eta_{t,h}^2)$, with $\eta_{t,h} = \sum_{j=0}^{h-1} b_j \varepsilon_{t+h-j}$, $b_j = \sum_{i=0}^{j} c_i$, $c_0 = 1$ and $c_j, j \geq 1$, satisfying $1 + \sum_{j=1}^{\infty} c_j z^j = 1/\alpha(z)$ (note that $\alpha(z)$ is defined after (1.1)). The first term on the right-hand sides of (1.4) and (1.5), originating from the random disturbances $\{\varepsilon_t\}$, is common for each multistep predictor, whereas the second terms on the right-hand sides of (1.4) and (1.5), arising from the estimation uncertainty, can vary with different $k$, different prediction methods and different parameter values. However, since only rates of convergence of the second terms are reported, (1.4) and (1.5) fail to depict these features, which are indispensable in performing predictor comparisons. To remedy this difficulty, the constants associated with the terms of order $n^{-1}$ in $E\{R_{P,n}(k)\}$ and $E\{R_{D,n}(k)\}$ need to be characterized. Recently, Ing [8] made a first step toward this goal. In the special case where $p = 0$ in (1.1) (the random walk model) and $k = h = 1$, he showed that

$$\lim_{n \to \infty} n(\text{MSPE}P_{n,1}(1) - \sigma^2) = \lim_{n \to \infty} E\left\{\frac{x_n^2}{n} n^2(\hat{\mathbf{a}}_n(1,1) - 1)^2\right\} = 2\sigma^2. \tag{1.6}$$

The main obstacle in dealing with the above expectation, as argued by Ing, is the fact that the square of the normalized regressor, $x_n^2/n$, and the square of the normalized estimator, $n^2(\hat{\mathbf{a}}_n(1,1) - 1)^2$, are **not** asymptotically independent – a situation somewhat different from that encountered in the stationary case. While Ing was able to overcome this difficulty, his approach, focusing only on the random walk model and the case of one-step-ahead prediction, cannot be directly applied to more general non-stationary AR models or multistep prediction cases.

Another subtle problem, related to the direct method, can be illustrated using the following special case of (1.1):

$$(1-B)(1 + 0.1B + 0.91B^2)x_{t+1} = (1 - 0.9B + 0.81B^2 - 0.91B^3)x_{t+1} = \varepsilon_{t+1}, \tag{1.7}$$

where $B$ is the back shift operator. Simple algebra yields

$$x_{t+1} = 0.181x_{t-2} + 0.819x_{t-3} + \varepsilon_{t+1} + 0.9\varepsilon_t. \tag{1.8}$$

As observed in (1.8), the direct method only requires two regressors to make a three-step-ahead prediction, which indicates the interesting fact that the minimal correct order for the direct method, determined by the prediction lead time and unknown parameters, can be strictly less than that for the plug-in method. In general, model (1.1) can be rewritten as

$$x_{t+h} = (A^{h-1}(p+1)\mathbf{a}(p+1))'\mathbf{x}_t(p+1) + \eta_{t,h}, \qquad h \geq 1,$$

where $\mathbf{a}(k) = (a_1, \ldots, a_k)'$, with $a_j = 0$ for $j > p + 1$,

$$A(k) = \left( \mathbf{a}(k) \,\middle|\, \begin{array}{c} I_{k-1} \\ \hline \mathbf{0}'_{k-1} \end{array} \right),$$

and $A^0(k) = I_k$. Let $\mathbf{a}(h, p+1) = (a_1(h, p+1), \ldots, a_{p+1}(h, p+1))' = A^{h-1}(p+1)\mathbf{a}(p+1)$. The above example leads us to define the minimal correct order for the $h$-step direct method, $p_h = \max\{j : 1 \leq j \leq p + 1, a_j(h, p+1) \neq 0\}$. As will be seen in Section 2 below, comparison results between the plug-in and direct predictors are very complicated in situations where $p_h < p_1$.

In Section 2, we first derive asymptotic expressions for MSPE$P_{n,h}(k_1)$ and MSPE$D_{n,h}(k_2)$ up to terms of order $n^{-1}$, where $k_1 \geq p_1$ and $k_2 \geq p_h$. The constants associated with the terms of order $n^{-1}$ in these expressions characterize how the prediction error is influenced by the orders, methods (plug-in or direct), values of parameters and even the unit roots. Based on these expressions, a series of examples (Examples 1–3) is given to illuminate that to find the asymptotically optimal (from the MSPE point of view) multistep predictor among candidate plug-in and direct predictors, prediction orders and prediction methods must simultaneously be taken into account. The traditional order selection criteria can no longer serve that purpose. Section 3 is devoted to alleviating this difficulty. Our strategy is to find a statistic for each MSPE$P_{n,h}(k)$ and MSPE$D_{n,h}(k), k = 1, \ldots, K$ and show that the ordering of these statistics coincides with the ordering of their corresponding multistep MSPEs. Here, $K \geq p_1$ is a known integer. In view of Ing [10], the statistics adopted in this section are the multistep generalizations of accumulated prediction errors (APEs) based on sequential plug-in and direct predictors, namely,

$$\mathrm{APE}P_{n,h}(k) = \sum_{i=m_h}^{n-h} (x_{i+h} - \hat{x}_{i+h}(k))^2 \tag{1.9}$$

and

$$\mathrm{APE}D_{n,h}(k) = \sum_{i=m_h}^{n-h} (x_{i+h} - \check{x}_{i+h}(k))^2, \tag{1.10}$$

where $m_h$ denotes the smallest positive number such that $\hat{\mathbf{a}}_i(h, K)$ and $\breve{\mathbf{a}}_i(h, K)$ are well defined for all $i \geq m_h$. Note that $\text{APE}P_{n,1}$ was first proposed by Rissanen [19]. A complete asymptotic analysis of $\text{APE}P_{n,1}$ was given by Wei [21, 22] under a model more general than (1.1). However, due to some "nice" properties in $\text{APE}P_{n,1}$ that are missing in its multistep counterparts (see Remarks 2 and 3 in Section 3), the asymptotic analysis of (1.9) and (1.10) in non-stationary AR processes is still lacking. We propose a resolution to this problem, which shows that every $\text{APE}P_{n,h}(k_1)$ and $\text{APE}D_{n,h}(k_2)$, with $k_1 \geq p_1$ and $k_2 \geq p_h$, can be asymptotically decomposed into two terms; one of which, due to estimation uncertainty, is of order $\log n$, and the other, due to the random disturbances, is of order $n$ and common for each predictor. More important, the constant associated with the term of $\log n$ in $\text{APE}P_{n,h}$ ($\text{APE}D_{n,h}$) is exactly the same as the one associated with the term of $n^{-1}$ in its corresponding $\text{MSPE}P_{n,h}$ ($\text{MSPE}D_{n,h}$). This special feature enables us to show that Ing's [10] asymptotically efficient predictor selection procedure (based on $\text{APE}P_{n,h}$ and $\text{APE}D_{n,h}$) in stationary AR processes can carry over to non-stationary cases and hence leads to a unified approach. Note that a predictor selection procedure is said to be asymptotically efficient if, with probability 1, it can choose the order/method combination with the minimal MSPE for all sufficiently large $n$; see Section 3 for the exact definition.

Despite its theoretical advantage, Ing's procedure suffers from unsatisfactory finite-sample performance, as explained at the beginning of Section 4. To fix this flaw, a new predictor selection method is proposed in Section 4. This new method not only shares the same asymptotic advantage as Ing's procedure, it also has satisfactory finite-sample performance, which is illustrated at the end of Section 4 through a simulation experiment. Appendices A–C contain the proofs of the theorems in Sections 2–4, respectively.

# 2. MSPEs of plug-in and direct predictors in the presence of unit roots

Throughout this section, it is assumed that in model (1.1) the $\varepsilon_t$'s are independent random variables with zero means and variances $\sigma^2 > 0$. Moreover, there are small positive numbers $\alpha_1$ and $\delta_1$ and a large positive number $M_1$ such that for $0 \leq s - \nu \leq \delta_1$

$$\sup_{1 \leq m \leq t < \infty, \|\mathbf{v}_m\| = 1} |F_{t,m,\mathbf{v}_m}(s) - F_{t,m,\mathbf{v}_m}(\nu)| \leq M_1(s - \nu)^{\alpha_1}, \tag{2.1}$$

where $\mathbf{v}_m = (v_1, \ldots, v_m)' \in R^m$, $\|\mathbf{v}_m\|^2 = \sum_{j=1}^m v_j^2$ and $F_{t,m,\mathbf{v}_m}(\cdot)$ denotes the distribution of $\sum_{l=1}^m v_l \varepsilon_{t+1-l}$.

In the case, where $\varepsilon_t$'s are i.i.d., the following lemma provides sufficient conditions under which (2.1) is fulfilled. The proof of this lemma can be found in Ing and Sin [12].

**Lemma 2.1.** *Let $\varepsilon_t$'s be i.i.d. random variables satisfying $E(\varepsilon_1) = 0, E(\varepsilon_1^2) > 0$, and $E(|\varepsilon_1|^\alpha) < \infty$ for some $\alpha > 2$. Assume also that for some positive constant $M_2 < \infty$,*

$$\int_{-\infty}^{\infty} |\varphi(t)| \, \mathrm{d}t \leq M_2, \tag{2.2}$$

*where $\varphi(t) = E\{\exp(\mathrm{i}t\varepsilon_1)\}$ is the characteristic function of $\varepsilon_1$. Then, for all $-\infty < t < \infty$, $m \geq 1$, $\mathbf{r}_m \in R^m$ and $\|\mathbf{r}_m\| = 1$, there is a finite positive constant $M_3$ such that*

$$\sup_{-\infty < x < \infty} f_{t,m,\mathbf{r}_m}(x) < M_3,$$

*where $f_{t,m,\mathbf{r}_m}(\cdot)$ is the density function of $(\varepsilon_t, \ldots, \varepsilon_{t+1-m})\mathbf{r}_m$. As a result, (2.1) follows.*

Since (2.2) is satisfied by most absolutely continuous distributions, (2.1) is flexible enough to accommodate a wide range time series applications. Note that (2.1) is given to ensure that the inverses of the normalized Fisher information matrices, $\hat{R}_n^{-1}(k)$ and $\bar{R}_{n,h}^{-1}(k)$, have finite positive moments in the senses of (A.1) and (A.19) (in Appendix A), where

$$\hat{R}_n(k) = \frac{1}{n} D_n(k) \sum_{j=k}^{n-1} \mathbf{x}_j(k) \mathbf{x}_j'(k) D_n(k)'$$

and

$$\bar{R}_{n,h}(k) = \frac{1}{n} \bar{D}_n(k) \sum_{j=k}^{n-h} \mathbf{x}_j(k) \mathbf{x}_j'(k) \bar{D}_n(k)',$$

with

$$D_n(k) = \begin{pmatrix} 1 & -1 & 0 & \cdots & 0 \\ 0 & 1 & -1 & \ddots & \vdots \\ \vdots & \ddots & \ddots & \ddots & 0 \\ 0 & \cdots & 0 & 1 & -1 \\ \frac{1}{\sqrt{n}} & \frac{-\alpha_1}{\sqrt{n}} & \cdots & \cdots & \frac{-\alpha_{k-1}}{\sqrt{n}} \end{pmatrix},$$

$\alpha_j = 0$ for $j > p$ and $\bar{D}_n(k)$ equal to $D_n(k)$ with $\alpha_i$ replaced by 0 for $i = 1, \ldots, k-1$. These results will be used to deal with the asymptotic properties of MSPE$P_{h,n}$ and MSPE$D_{h,n}$; see the proofs of Theorems 2.2 and 2.3 for details. Theorems 2.2 below provides an asymptotic expression for MSPE$P_{n,h}(k)$ with $k \geq p_1$. Before stating the result, we need to define $S_M^0(k) = I_k$ and with $\alpha(k) = (\alpha_1, \ldots, \alpha_k)'$,

$$S_M(k) = \left( \alpha(k) \,\middle|\, \begin{array}{c} I_{k-1} \\ \hline \mathbf{0}_{k-1}' \end{array} \right).$$

**Theorem 2.2.** *Assume that $\{x_t\}$ satisfies model (1.1). Also assume that $\{\varepsilon_t\}$ satisfies (2.1) and*

$$E(|\varepsilon_1|^{\theta_h}) < \infty,$$

*where $\theta_h = \max\{8, 2(h+2)\} + \delta$ for some $\delta > 0$. Then, for $k \geq p_1$ and $h \geq 1$,*

$$n(\mathrm{MSPE}P_{n,h}(k) - \sigma_h^2) = 2\sigma^2 \left(\sum_{j=0}^{h-1} b_j\right)^2 + f_{1,h}(k-1) + \mathrm{o}(1), \tag{2.3}$$

*where $f_{1,h}(0) = 0$ and for $k \geq 2$,*

$$f_{1,h}(k-1) = \mathrm{tr}(\Gamma(k-1)M_h(k-1)\Gamma^{-1}(k-1)M_h'(k-1))\sigma^2,$$

*with $M_h(k-1) = \sum_{j=0}^{h-1} b_j S_M^{h-1-j}(k-1)$, $\Gamma(k-1) = \lim_{j \to \infty} E(\mathbf{s}_j(k-1)\mathbf{s}_j'(k-1))$, $\mathbf{s}_j(k-1) = (s_j, \ldots, s_{j-k+2})'$ and $s_j = x_j - x_{j-1}$.*

An asymptotic expression for $\mathrm{MSPE}D_{n,h}(k)$, with $k \geq p_h$, is given as follows:

**Theorem 2.3.** *Let the assumptions of Theorem 2.2 hold, with $\theta_h$ replaced by $8 + \delta$ for some $\delta > 0$. Then, for $k \geq p_h$ and $h \geq 1$,*

$$n(\mathrm{MSPE}D_{n,h}(k) - \sigma_h^2) = 2\sigma^2 \left(\sum_{j=0}^{h-1} b_j\right)^2 + f_{2,h}(k-1) + \mathrm{o}(1), \tag{2.4}$$

*where $f_{2,h}(0) = 0$, for $k \geq 2$,*

$$f_{2,h}(k-1) = \mathrm{tr}\left\{\Gamma^{-1}(k-1) \lim_{t \to \infty} \mathrm{cov}\left(\sum_{j=0}^{h-1} b_j \mathbf{s}_{t+j}(k-1)\right)\right\}\sigma^2,$$

*and for random vector $\mathbf{y}$, $\mathrm{cov}(\mathbf{y}) = E\{(\mathbf{y} - E(\mathbf{y}))(\mathbf{y} - E(\mathbf{y}))'\}$.*

Theorems 2.2 and 2.3 show that each $n(\mathrm{MSPE}P_{n,h}(k_1) - \sigma_h^2)$ and $n(\mathrm{MSPE}D_{n,h}(k_2) - \sigma_h^2)$, with $k_1 \geq p_1$ and $k_2 \geq p_h$, can be asymptotically decomposed as a sum of two terms. The first term, $2\sigma^2(\sum_{j=0}^{h-1} b_j)^2$, arising from predicting the non-stationary component in model (1.1), is common for each predictor, whereas the second term, $f_{1,h}(k-1)$ (or $f_{2,h}(k-1)$), arising from predicting the stationary component in model (1.1), can vary with different orders and methods. The following examples help provide a better understanding of Theorems 2.2 and 2.3.

**Example 1.** When $k \geq \max\{2, p_1\}$ and $h = 2$, by (2.3) and (2.4), it is straightforward to show that

$$f_{1,2}(k-1) = \{(k-2) + \alpha_{k-1}^2 + 2\alpha_1 b_1 + b_1^2(k-1)\}\sigma^2 \tag{2.5}$$

and

$$f_{2,2}(k-1) = \{(k-1)(1 + b_1^2) + 2\alpha_1 b_1\}\sigma^2, \tag{2.6}$$

which yields

$$f_{2,2}(k-1) - f_{1,2}(k-1) = (1 - \alpha_{k-1}^2)\sigma^2 > 0. \tag{2.7}$$

Moreover, by an argument similar to that used to prove (17) of Ing [9], it can be shown that for $k \geq \max\{2, p_1\}$ and $h \geq 2$,

$$f_{2,h}(k-1) - f_{1,h}(k-1) \geq f_{2,2}(k-1) - f_{1,2}(k-1) > 0, \tag{2.8}$$

and hence $\hat{x}_{n+h}(k)$ is asymptotically more efficient than $\check{x}_{n+h}(k)$ in this case.

As shown in Section 1, it is possible that $p_h < p_1$. In this case, it would be more interesting to compare $n(\text{MSPE}P_{n,h}(p_1) - \sigma_h^2)$ and $n(\text{MSPE}D_{n,h}(p_h) - \sigma_h^2)$ rather than those MSPEs of the same order. The following example shows that the advantage of the plug-in predictor illustrated in Example 1 vanishes in this kind of comparison.

**Example 2.** Assume

$$(1 - B)(1 + a_1 B + \cdots + a_p B^p)x_t = e_t,$$

where $p \geq 2$, $1 + a_1 z + \cdots + a_p z^p \neq 0$ for $|z| \leq 1$ and $a_p \neq 0$. If $a_1 = 1$, then it is not difficult to see that $p_2 = p_1 - 1 = p$ and $f_{2,2}(p) - f_{2,2}(p-1) = \sigma^2$. In addition, (2.7) implies $f_{2,2}(p) - f_{1,2}(p) = (1 - a_p^2)\sigma^2$. As a result,

$$n\{\text{MSPE}P_{n,2}(p_1) - \sigma^2\} - n\{\text{MSPE}D_{n,2}(p_2) - \sigma^2\} \to a_p^2\sigma^2 > 0,$$

as $n \to \infty$. Hence $\check{x}_{n+2}(p_2)$ is asymptotically more efficient than $\hat{x}_{n+2}(p_1)$ in this case.

When $h = 2$ and $p_1 \geq 2$, Examples 1 and 2 together suggest a simple rule that $\hat{x}_{n+2}(p_1)$ is asymptotically more efficient than $\check{x}_{n+2}(p_2)$ if $p_1 = p_2$; and the conclusion is reversed if $p_1 > p_2$. This rule, however, fails to hold for $h \geq 3$, as detailed in the following example.

**Example 3.** Consider the following AR(4) model

$$(1 - B)(1 + a_1 B)(1 + a_2 B^2)x_t$$
$$= \{1 - (1 - a_1)B - (a_1 - a_2)B^2 - a_2(1 - a_1)B^3 - a_1 a_2 B^4\}x_t = e_t,$$

**Table 1.** The values of Diff $= f_{2,3}(2) - f_{1,3}(3)$

| $a_1$ | 0.1 | 0.2 | 0.3 | 0.4 | 0.5 | 0.6 | 0.7 | 0.8 | 0.9 |
|-------|-----|-----|-----|-----|-----|-----|-----|-----|-----|
| Diff | $-0.378$ | $-0.013$ | $0.197$ | $0.310$ | $0.354$ | $0.336$ | $0.247$ | $0.051$ | $-0.321$ |

where $0 < a_1 < 1$ and $a_2 = a_1^2 - a_1 + 1$. It is straightforward to show that $p_3 = 3 = p_1 - 1$. By numerical calculations, we obtain the values of $f_{2,3}(2) - f_{1,3}(3)$, with $a_1 = 0.1, 0.2, \ldots, 0.9$; see Table 1. According to Table 1, $\check{x}_{n+3}(p_3)$ is asymptotically more efficient than $\hat{x}_{n+3}(p_1)$ in cases of $a_1 = 0.1, 0.2, 0.9$, and less efficient than $\hat{x}_{n+3}(p_1)$ in all other cases.

Consequently, when $h \geq 3$, the rankings of $\hat{x}_{n+h}(p_1)$ and $\check{x}_{n+h}(p_h)$ are determined not only by whether $p_h < p_1$, but also by the values of the unknown parameters. Simply determining $p_1$ or $p_h$ through certain consistent model selection techniques cannot guarantee optimal multistep prediction (from the MSPE point of view) in situations, where plug-in and direct predictors are simultaneously taken into account. This phenomenon was first reported by Ing [10] in stationary AR models. The above three examples show that the same difficulty occurs in the presence of unit root. In the next two sections, some proposals toward resolving this problem are given.

## 3. Multistep accumulated prediction errors

Let $\hat{x}_{n+h}(k), k = 1, \ldots, K$ and $\check{x}_{n+h}(k), k = 1, \ldots, K$, be candidate plug-in and direct predictors, where $h \geq 1$ and $K \geq p_1$. For convenience, we use $(k, 1)$ to denote $\hat{x}_{n+h}(k)$ and $(k, 2)$ to denote $\check{x}_{n+h}(k)$. In response to the difficulty mentioned at the end of the previous section, this section attempts to choose the order/method combination having the minimal MSPE instead of identifying $p_1$ or $p_h$. To this end, the loss functions of $(k, 1)$ and $(k, 2)$ are defined to be

$$L_{1,h}(k) = \begin{cases} \lim_{n \to \infty} n(\text{MSPE}P_{n,h}(k) - \sigma_h^2), & \text{if } p_1 \leq k \leq K, \\ \infty, & \text{if } k < p_1, \end{cases} \tag{3.1}$$

and

$$L_{2,h}(k) = \begin{cases} \lim_{n \to \infty} n(\text{MSPE}D_{n,h}(k) - \sigma_h^2), & \text{if } p_h \leq k \leq K, \\ \infty, & \text{if } \quad k < p_h, \end{cases} \tag{3.2}$$

respectively. Note that the existence of the above limits is ensured by Theorems 2.2 and 2.3; and in order to have the prediction loss due to underspecification be much larger than the one due to overspecification, the loss function values of $(k, 1)$ with $k < p_1$ and $(k, 2)$ with $k < p_h$ are set to $\infty$. A predictor selection criterion, $(\tilde{k}_n, \tilde{j}_n)$, with $1 \leq \tilde{k}_n \leq K$ and $1 \leq \tilde{j}_n \leq 2$, is said to be asymptotically efficient if

$$P((\tilde{k}_n, \tilde{j}_n) \in C_{h,K}, \text{eventually}) = 1, \tag{3.3}$$

where

$$C_{h,K} = \left\{ (k, j) : 1 \leq k \leq K, 1 \leq j \leq 2 \text{ and } L_{j,h}(k) = \min_{1 \leq k_0 \leq K, 1 \leq j_0 \leq 2} L_{j_0,h}(k_0) \right\}.$$

Therefore, with probability 1 $(\tilde{k}_n, \tilde{j}_n)$ can choose the predictor having the minimal loss function value for all sufficiently large $n$.

The goal of this section is to show that (3.3) is fulfilled by $(\hat{k}_n, \hat{j}_n)$. Here, $(\hat{k}_n, \hat{j}_n)$, first proposed by Ing [10], is obtained through the following procedure:

**Step 1**. Define $\hat{k}_{D,n}^{(1)} = \arg\min_{1 \le k \le K} \mathrm{APE}D_{n,1}(k)$.

**Step 2**. Define

$$\hat{k}_{D,n}^{(h)} = \arg\min_{1 \le k \le K} \mathrm{APE}D_{n,h}(k)$$

and define

$$\hat{k}_n^{(1,h)} = \arg\min_{\hat{k}_{D,n}^{(1)} \le k \le K} \mathrm{APE}P_{n,h}(k).$$

**Step 3**. If $\mathrm{APE}D_{n,h}(\hat{k}_{D,n}^{(h)}) > \mathrm{APE}P_{n,h}(\hat{k}_n^{(1,h)})$, then $(\hat{k}_n, \hat{j}_n) = (\hat{k}_n^{(1,h)}, 1)$; otherwise $(\hat{k}_n, \hat{j}_n) = (\hat{k}_{D,n}^{(h)}, 2)$.

*Remark 1.* Our analysis below implies that the asymptotic properties of $(\hat{k}_n, \hat{j}_n)$ remain unchanged if Step 1 is skipped and $\hat{k}_n^{(1,h)}$ in Step 2 is defined to be $\arg\min_{1 \le k \le K} \mathrm{APE}P_{n,h}(k)$.

In the sequel, the above procedure will be referred to as Procedure I. We begin by investigating the asymptotic properties of $\mathrm{APE}P_{n,h}(k)$ and $\mathrm{APE}D_{n,h}(k)$ in the correctly specified case. Note that for $k \ge p_1$,

$$\mathrm{APE}P_{n,h}(k) = \sum_{i=m_h}^{n-h} \{\eta_{i,h} - \mathbf{x}_i'(k)\hat{L}_{i,h}(k)(\hat{\mathbf{a}}_i(1,k) - \mathbf{a}(k))\}^2, \tag{3.4}$$

where $\hat{L}_{i,h}(k) = \sum_{j=0}^{h-1} b_j \hat{A}_i^{h-1-j}(k)$; and for $k \ge p_h$,

$$\mathrm{APE}D_{n,h}(k) = \sum_{i=m_h}^{n-h} \{\eta_{i,h} - \mathbf{x}_i'(k)(\hat{\mathbf{a}}_i(h,k) - \mathbf{a}_D(h,k))\}^2, \tag{3.5}$$

where $\mathbf{a}_D(h,k) = (a_1(h,p+1), \ldots, a_k(h,p+1))'$, with $a_j(h,p+1), 1 \le j \le p+1$, defined in Section 1 and $a_j(h,p+1) = 0$ if $j > p+1$.

**Theorem 3.1.** *Assume that $\{x_t\}$ satisfies model (1.1) and $\{\varepsilon_t\}$ is a sequence of independent random noises with zero means and common variance $\sigma^2 > 0$. Moreover, assume $\sup_t E(|\varepsilon_t|^\alpha) < \infty$ for some $\alpha > 2$. Then, for $k \ge p_1$ and $h \ge 1$,*

$$\mathrm{APE}P_{n,h}(k) - \sum_{i=m_h}^{n-h} \eta_{i,h}^2 = \left\{ 2\sigma^2 \left( \sum_{j=0}^{h-1} b_j \right)^2 + f_{1,h}(k-1) \right\} \log n + \mathrm{o}(\log n) \qquad \text{a.s.}$$

$$= L_{1,h}(k) \log n + \mathrm{o}(\log n) \qquad \text{a.s.} \tag{3.6}$$

**Remark 2.** As shown in (B.18),

$$\text{APEP}_{n,h}(k) - \sum_{i=m_h}^{n-h} (\eta_{i,h})^2 = \sum_{i=m_h}^{n-h} \{\mathbf{x}_i'(k)\hat{L}_{i,h}(k)(\hat{\mathbf{a}}_i(1,k) - \mathbf{a}(k))\}^2(1+\text{o}(1))$$
$$+ \text{O}(1) \qquad \text{a.s.}$$

Therefore, the main task of proving (3.6) is to explore the almost sure properties of the first term on the right-hand side of the above equality. Through a recursive expression for $Q_n(1,k)$, where, with $V_i^{-1}(k) = \sum_{i=k}^{i} \mathbf{x}_i(k)\mathbf{x}_i'(k)$,

$$Q_n(1,k) = \sum_{i=k}^{n-1} \{\mathbf{x}_i'(k)(\hat{\mathbf{a}}_n(1,k) - \mathbf{a}(k))\}^2$$

$$= \left(\sum_{i=k}^{n-1} \mathbf{x}_i'(k)\varepsilon_{i+1}\right) V_{n-1}(k) \left(\sum_{i=k}^{n-1} \mathbf{x}_i(k)\varepsilon_{i+1}\right)$$

is the (second-order) residual sum of squares for one-step predictions, Lai and Wei [14] established a connection between $Q_n(1,k)$ and its sequential counterpart,

$$\sum_{i=m_1}^{n-1} \{\mathbf{x}_i'(k)(\hat{\mathbf{a}}_i(1,k) - \mathbf{a}(k))\}^2 = \sum_{i=m_h}^{n-1} \left\{\mathbf{x}_i'(k)V_{i-1}(k)\left(\sum_{j=k}^{i-1} \mathbf{x}_j(k)\varepsilon_{j+1}\right)\right\}^2. \quad (3.7)$$

Based on this connection and some strong laws for martingales, Wei [21, 22] subsequently obtained an asymptotic expression for the left-hand side of (3.6) in the case of $h = 1$. However, it is extremely difficult to obtain an analyzable recursive formula for the multistep analog of $Q_n(1,k)$, $Q_n(h,k) = \sum_{i=k}^{n-h} \{\mathbf{x}_i'(k)\hat{L}_{n,h}(k)(\hat{\mathbf{a}}_n(1,k) - \mathbf{a}(k))\}^2, h \geq 2$, due to the appearance of $\hat{L}_{n,h}(k)$. Hence, Wei's approach is not easily extended to the case of multistep predictions. By observing

$$Q_n(h,k) = \left(\sum_{i=k}^{n-1} \mathbf{x}_i'(k)\varepsilon_{i+1}\right) S_n'(k)V_{n-h}(k)S_n(k)\left(\sum_{i=k}^{n-1} \mathbf{x}_i(k)\varepsilon_{i+1}\right),$$

where

$$S_n(k) = \left(\sum_{i=k}^{n-h} \mathbf{x}_i(k)\mathbf{x}_i'(k)\right) \hat{L}_{n,h}(k) \left(\sum_{i=k}^{n-1} \mathbf{x}_i(k)\mathbf{x}_i'(k)\right)^{-1},$$

Ing [10], under stationary AR processes, adopted

$$Q_n^*(h,k) = \left(\sum_{i=k}^{n-1} \mathbf{x}_i'(k)\varepsilon_{i+1}\right) S'(k)V_{n-h}(k)S(k)\left(\sum_{i=k}^{n-1} \mathbf{x}_i'(k)\varepsilon_{i+1}\right)$$

to replace $Q_n(h, k)$, where $S(k)$ is the almost sure limit of $S_n(k)$ that is a non-random matrix. He then obtained a recursive formula for $Q_n^*(h, k)$ and established a connection between $Q_n^*(h, k)$ and $\sum_{i=k}^{n-h} \{\mathbf{x}_i'(k) L_{i,h}(k)(\hat{\mathbf{a}}_i(1, k) - \mathbf{a}(k))\}^2$, which further yields an asymptotic expression for the latter. Unfortunately, when model (1.1) is assumed, $S_n(k)$, with $k \geq 2$, no longer has an almost sure and non-random limit, which makes it hard to apply Ing's [10] approach to the non-stationary case. To obtain (3.6), extra effort is made to overcome the above difficulties; see Appendix B for details. For some other interesting analysis of APEs in various non-standard situations, see de Luna and Skouras [4] and Bercu [2].

**Theorem 3.2.** *Let the assumptions of Theorem 3.1 hold. Then, for $k \geq p_h$ and $h \geq 1$,*

$$\text{APE}D_{n,h}(k) - \sum_{i=m_h}^{n-h} \eta_{i,h}^2 = \left\{ 2\sigma^2 \left( \sum_{j=0}^{h-1} b_j \right)^2 + f_{2,h}(k-1) \right\} \log n + \mathrm{o}(\log n) \qquad a.s.$$

$$= L_{2,h}(k) \log n + \mathrm{o}(\log n) \qquad a.s.$$
(3.8)

**Remark 3.** As indicated in (B.35),

$$\text{APE}D_{n,h}(k) - \sum_{i=m_h}^{n-h} \eta_{i,h}^2 = (1 + \mathrm{o}(1)) \sum_{i=m_h}^{n-h} \left\{ \mathbf{x}_i'(k) V_{i-h}(k) \left( \sum_{j=k}^{i-h} \mathbf{x}_j(k) \eta_{j,h} \right) \right\}^2$$

$$+ \mathrm{O}(1) \qquad a.s.$$

While

$$\sum_{i=m_h}^{n-h} \left\{ \mathbf{x}_i'(k) V_{i-h}(k) \left( \sum_{j=k}^{i-h} \mathbf{x}_j(k) \eta_{j,h} \right) \right\}^2$$

looks very similar to (3.7), Wei's approach for the one-step APE still cannot be applied to it because $\sum_{j=k}^{i-h} \mathbf{x}_j(k) \eta_{j,h}, h \geq 2$, is not a martingale transformation. While Ing [10] resolved this difficulty in the stationary AR model, his method, which is highly reliant on the stationary assumption, is not applicable to the unit root processes.

**Remark 4.** Theorems 2.2, 2.3, 3.1 and 3.2 together disclose a fascinating fact that the constants associated with the terms of order $n^{-1}$ in $\text{MSPE}P_{n,h}(k_1)$ and $\text{MSPE}D_{n,h}(k_2)$, with $k_1 \geq p_1$ and $k_2 \geq p_h$, are exactly the same as the constants associated with the terms of order $\log n$ in their corresponding multistep APEs. While $\text{MSPE}P_{n,h}(k_1)$ and $\text{MSPE}D_{n,h}(k_2)$ are unobservable, this special property allows us to preserve their asymptotic rankings through the values of the associated multistep APEs, which can be easily obtained from the data. This is also the driving motivation for constructing $(\hat{k}_n, \hat{j}_n)$ in model (1.1).

Before showing the asymptotic efficiency of $(\hat{k}_n, \hat{j}_n)$, we need to investigate the asymptotic properties of $\text{APED}_{n,k}(k)$ in misspecified cases.

**Theorem 3.3.** *Let the assumptions of Theorem 3.1 hold. Then, for $1 \le k < p_h$ and $h \ge 1$,*

$$\liminf_{n \to \infty} \frac{1}{n} \left( \text{APED}_{n,h}(k) - \sum_{j=m_h}^{n-h} \eta_{j,h}^2 \right) > 0 \qquad \text{a.s.} \tag{3.9}$$

We are now in a position to state the main result of this section.

**Theorem 3.4.** *Let the assumptions of Theorem 3.1 hold. Then, for $K \ge p_1$, $(\hat{k}_n, \hat{j}_n)$ is asymptotically efficient in the sense of (3.3).*

**Remark 5.** Since Ing [10] showed that $(\hat{k}_n, \hat{j}_n)$ is also asymptotically efficient in stationary AR models, Theorem 3.4, together with Ing's result, provides a unified approach for choosing the (asymptotically) optimal multistep predictor for AR processes with or without unit roots. While it is possible to select multistep predictors after unit root tests are performed (which means that the selection procedure will be carried out based on the differenced data if the unit-root hypothesis is not rejected), all unit root tests suffer from low power when the process is near unity. One can hardly expect a reliable selection/prediction result once the process is erroneously differenced.

Before leaving this section, we note that to analyze the effect of the estimation of the mean into the performance of the predictors, one may consider a unit root AR model with drift,

$$A(B)x_{t+1} = \beta + \varepsilon_{t+1}, \tag{3.10}$$

where $A(B)$ is defined after (1.1) and $-\infty < \beta < \infty$ is some real number. In the case of $h = 1$, we have obtained (through non-trivial modifications of the proofs of the results in Sections 2 and 3) that if $\beta \ne 0$, then for $k \ge p_1$,

$$\lim_{n \to \infty} n\{E(x_{n+1} - \hat{x}_{n+1}(k))^2 - \sigma^2\} = (k+3)\sigma^2,$$

and

$$\text{APEP}_{n,1}(k) - \sum_{i=m_1}^{n-1} \varepsilon_{i+1}^2 = \sigma^2(k+3)\log n + \text{o}(\log n) \qquad \text{a.s.},$$

where $\hat{x}_{n+1}(k) = \check{x}_{n+1}(k) = \mathbf{w}'_n(k)\hat{\mathbf{a}}_n(1, k)$, with $\mathbf{w}_j(k) = (1, \mathbf{x}'_j(k))'$ and $\hat{\mathbf{a}}_j(1, k)$ satisfying $\{\sum_{l=k}^{j-1} \mathbf{w}_l(k)\mathbf{w}'_l(k)\}\hat{\mathbf{a}}_j(1, k) = \sum_{l=k}^{j-1} \mathbf{w}_l(k)x_{l+1}$. Moreover, if $\beta = 0$, then for $k \ge p_1$,

$$\lim_{n \to \infty} n\{E(x_{n+1} - \hat{x}_{n+1}(k))^2 - \sigma^2\} = (k+2)\sigma^2,$$

and

$$\mathrm{APE}P_{n,1}(k) - \sum_{i=m_1}^{n-1} \hat{\varepsilon}_{i+1}^2 = \sigma^2(k+2)\log n + \mathrm{o}(\log n) \qquad \text{a.s.}$$

As observed in the above four equalities, the correspondence between APE and MSPE remains valid under model (3.10), regardless of whether $\beta = 0$ or not. Therefore, it is natural to conjecture that under model (3.10), (i) this correspondence can be extended to the case of $h > 1$; and (ii) Procedure I is still asymptotically efficient for multistep prediction. However, we shall not pursue a proof of these conjectures here, since it goes beyond the scope of this paper.

# 4. New criteria

Although Theorem 3.4 shows that $(\hat{k}_n, \hat{j}_n)$ is asymptotically efficient in the sense of (3.3), surprisingly, its finite sample performance is rather unsatisfactory. Simulation results show that the rankings of $\mathrm{APE}P_{n,h}(k_1)$ and $\mathrm{APE}D_{n,h}(k_2)$ are often inconsistent with the rankings $L_{1,h}(k_1)$ and $L_{2,h}(k_2)$ even when $n > 500$. One possible explanation for this phenomenon is as follows: In view of (3.4) and (3.5), for $k_1 \geq p_1$ and $k_2 \geq p_h$,

$$\mathrm{APE}P_{n,h}(k_1) - \mathrm{APE}D_{n,h}(k_2)$$

$$= \sum_{i=m_h}^{n-h} \{\mathbf{x}_i'(k_1)\hat{L}_{i,h}(k_1)(\hat{\mathbf{a}}_i(1,k_1) - \mathbf{a}(k_1))\}^2$$

$$- \sum_{i=m_h}^{n-h} \{\mathbf{x}_i'(k_2)(\bar{\mathbf{a}}_i(h,k_2) - \mathbf{a}_D(h,k_2))\}^2 - 2\sum_{i=m_h}^{n-h} \mathbf{x}_i'(k_1)\hat{L}_{i,h}(k_1)(\hat{\mathbf{a}}_i(1,k_1) - \mathbf{a}(k_1))\eta_{i,h}$$

$$+ 2\sum_{i=m_h}^{n-h} \mathbf{x}_i'(k_2)(\bar{\mathbf{a}}_i(h,k_2) - \mathbf{a}_D(h,k_2))\eta_{i,h} \equiv (I) - (II) - (III) + (IV).$$
(4.1)

While the cross-product terms, $(III)$ and $(IV)$, in (4.1) are almost surely of order $\mathrm{o}(\log n)$ and asymptotically negligible compared to $(I)$ and $(II)$ (see Appendix B), we have found that the finite sample values of $(III)$ and $(IV)$ can differ remarkably. This "nonuniformity" feature causes "rank-distortion" when we perform cross-method comparisons.

To overcome the above difficulty, we consider using $\mathrm{PMIC}_{n,h}(k)$ and $\mathrm{DMIC}_{n,h}(k)$ to replace $\mathrm{APE}P_{n,h}(k)$ and $\mathrm{APE}D_{n,h}(k)$ in Procedure I, where

$$\mathrm{PMIC}_{n,h}(k)$$

$$= \hat{\sigma}_{P,n}^2(h,k) + \mathrm{tr}\left\{\left(\sum_{j=k}^{n-h} \mathbf{x}_j(k)\mathbf{x}_j'(k)\right)\ddot{L}_{h,n}(k)\left(\sum_{j=k}^{n-h} \mathbf{x}_j(k)\mathbf{x}_j'(k)\right)^{-1}\ddot{L}_{h,n}'(k)\right\}\hat{\sigma}_n^2 C_n,$$
(4.2)

and

$$\text{DMIC}_{n,h}(k) = \hat{\sigma}_D^2(h,k) + \text{tr}\left\{\left(\sum_{j=k}^{n-h} \mathbf{x}_j(k)\mathbf{x}_j'(k)\right)^{-1}\left(\sum_{j=k}^{n-2h+1} \mathbf{z}_j(k)\mathbf{z}_j'(k)\right)\right\}\tilde{\sigma}_n^2 C_n, \quad (4.3)$$

where $\lim_{n\to\infty} C_n = 0$ and $\liminf_{n\to\infty} C_n n/\log n > 0$. Note that $\hat{\sigma}_{P,n}^2(h,k) = (n-h-K)^{-1}\sum_{j=K}^{n-h}\{x_{j+h} - \hat{\mathbf{a}}_n(h,k)\mathbf{x}_j(k)\}^2$ and $\hat{\sigma}_{D,n}^2(h,k) = (n-h-K)^{-1}\sum_{j=K}^{n-h}\{x_{j+h} - \tilde{\mathbf{a}}_n(h,k)\mathbf{x}_j(k)\}^2$ are the $h$-step residual mean squared errors obtained from the $k$-regressor plug-in and direct methods, respectively; $\tilde{\sigma}_n^2 = \hat{\sigma}_{P,n}^2(1,K) = \hat{\sigma}_{D,n}^2(1,K)$ is the one-step residual mean squared error obtained from the largest candidate model, $\mathbf{z}_j(k) = \sum_{i=0}^{h-1} \hat{b}_{i,n}\mathbf{x}_{j+i}(k)$, and $\ddot{L}_{h,n}(k) = \sum_{j=0}^{h-1} \hat{b}_{j,n}\hat{A}_n^{h-1-j}(k)$, where $\hat{b}_{0,n} = 1$, and for $j \geq 1$, $\hat{b}_{j,n} = \sum_{l=1}^{j} \hat{b}_{j-l,n}\hat{a}_{l,n}(1,K)$, with $(\hat{a}_{1,n}(1,K),\ldots,\hat{a}_{K,n}(1,K))' = \hat{\mathbf{a}}_n(1,K)$ and $\hat{a}_{l,n}(1,K) = 0$ if $l > K$.

Here, we briefly describe some of the theoretical rationale behind this new criterion. Observe that

$$\begin{aligned}
&\text{PMIC}_{n,h}(k_1) - \text{DMIC}_{n,h}(k_2) \\
&= \hat{\sigma}_{P,n}^2(h,k_1) - \hat{\sigma}_D^2(h,k_2) \\
&\quad + \text{tr}\left\{\left(\sum_{j=k_1}^{n-h} \mathbf{x}_j(k_1)\mathbf{x}_j'(k_1)\right)\ddot{L}_{h,n}(k_1)\left(\sum_{j=k_1}^{n-h} \mathbf{x}_j(k_1)\mathbf{x}_j'(k_1)\right)^{-1}\ddot{L}_{h,n}'(k_1)\right\}\tilde{\sigma}_n^2 C_n \\
&\quad - \text{tr}\left\{\left(\sum_{j=k_2}^{n-h} \mathbf{x}_j(k_2)\mathbf{x}_j'(k_2)\right)^{-1}\left(\sum_{j=k_2}^{n-2h+1} \mathbf{z}_j(k_2)\mathbf{z}_j'(k_2)\right)\right\}\tilde{\sigma}_n^2 C_n.
\end{aligned} \quad (4.4)$$

It is shown in Appendix C that when $k_1 \geq p_1$ and $k_2 \geq p_h$,

$$\begin{aligned}
&\text{tr}\left\{\left(\sum_{j=k_1}^{n-h} \mathbf{x}_j(k_1)\mathbf{x}_j'(k_1)\right)\ddot{L}_{h,n}(k_1)\left(\sum_{j=k_1}^{n-h} \mathbf{x}_j(k_1)\mathbf{x}_j'(k_1)\right)^{-1}\ddot{L}_{h,n}'(k_1)\right\}\tilde{\sigma}_n^2 \\
&\quad - \text{tr}\left\{\left(\sum_{j=k_2}^{n-h} \mathbf{x}_j(k_2)\mathbf{x}_j'(k_2)\right)^{-1}\left(\sum_{j=k_2}^{n-2h+1} \mathbf{z}_j(k_2)\mathbf{z}_j'(k_2)\right)\right\}\tilde{\sigma}_n^2 \\
&= L_{1,h}(k_1) - L_{2,h}(k_2) + \text{o}(1) \qquad \text{a.s.}
\end{aligned} \quad (4.5)$$

Therefore, the trace terms in (4.2) and (4.3) play roles in keeping the rankings of their corresponding loss functions. On the other hand, for $k_1 \geq p_1$ and $k_2 \geq p_h$, the weight associated with the trace terms, $C_n$, asymptotically dominates $\hat{\sigma}_{P,n}^2(h,k_1) - \hat{\sigma}_D^2(h,k_2)$ (see (C.2)), which helps to protect the trace term effects in (4.4) from being distorted by $\hat{\sigma}_{P,n}^2(h,k_1) - \hat{\sigma}_D^2(h,k_2)$. In fact, our simulations reveal that this domination usually occurs quite early (particular when $C_n$ is relatively large), and hence considerably alleviate the dilemma encountered by Procedure I in finite samples. (Note that $\hat{\sigma}_{P,n}^2(h,k)$ and $\hat{\sigma}^2_{D,n}(h,k)$ cannot be dropped from (4.2) and (4.3) because they are necessary for preventing underspecification; see, e.g., (C.1).) The following is the new predictor selection procedure (which is referred to as Procedure II) and its asymptotic property.

**Step 1**. Define $\hat{O}_n^{(1)} = \arg\min_{1 \leq k \leq K} \text{DMIC}_{n,1}(k)$.

**Step 2**. Define

$$\hat{O}_n^{(h)} = \arg\min_{1 \leq k \leq K} \text{DMIC}_{n,h}(k)$$

and define

$$\hat{O}_n^{(1,h)} = \arg\min_{\hat{O}_n^{(1)} \leq k \leq K} \text{PMIC}_{n,h}(k).$$

**Step 3**. If $\text{DMIC}_{n,h}(\hat{O}_n^{(h)}) > \text{PMIC}_{n,h}(\hat{O}_n^{(1,h)})$, then $(\hat{O}_n, \hat{M}_n) = (\hat{O}_n^{(1,h)}, 1)$; otherwise $(\hat{O}_n, \hat{M}_n) = (\hat{O}_n^{(h)}, 2)$.

**Theorem 4.1.** *Let the assumptions of Theorem 3.1 hold. Then, for $K \geq p_1$, $(\hat{O}_n, \hat{M}_n)$ is asymptotically efficient in the sense of (3.3).*

***Remark 6.*** Although (4.5) holds, it is worth mentioning that the trace terms in (4.5) are not consistent estimators of their corresponding loss functions $L_{1,h}(k_1)$ and $L_{2,h}(k_2)$; see (C.5) and (C.6) in Appendix C.

***Remark 7.*** Following an argument similar to that used in the proof of Theorem 4.1, it is not difficult to show that $(\hat{O}_n, \hat{M}_n)$ is also asymptotically efficient in stationary AR models.

To illustrate the asymptotic results obtained in Theorem 4.1, we conduct a simulation study. The data generating processes (DGPs) are given by

**DGP I** $x_t = -0.8x_{t-2} + \varepsilon_t$,
**DGP II** $x_t = 0.3x_{t-1} - 0.8x_{t-2} + \varepsilon_t$,
**DGP III** $x_t = 0.2x_{t-2} + 0.8x_{t-3} + \varepsilon_t$,
**DGP IV** $x_t = 0.3x_{t-1} - 0.1x_{t-2} + 0.8x_{t-3} + \varepsilon_t$,
**DGP V** $x_t = 0.9x_{t-1} - 0.81x_{t-2} + \varepsilon_t$,
**DGP VI** $x_t = 0.6x_{t-1} - 0.36x_{t-2} + \varepsilon_t$,
**DGP VII** $x_t = 0.9x_{t-1} - 0.81x_{t-2} + 0.91x_{t-3} + \varepsilon_t$,
**DGP VIII** $x_t = 0.9x_{t-1} - 0.56x_{t-2} + 0.66x_{t-3} + \varepsilon_t$,

where $\varepsilon_t$'s are independent and identically $\mathcal{N}(0, 25)$ distributed. We aim to select two-step ($h = 2$) predictors for DGPs I–IV and three-step ($h = 3$) predictors for DGPs V–VIII using Procedure II with $C_n = \log n/n$, $2\log n/n$ and $3\log n/n$, which will be referred to as Procedures A, B and C, respectively. The candidate predictors are set to $(i, j)$, $i = 1, \ldots, 10$ and $j = 1, 2$. According to Section 2 and Section 2 of Ing [10], the asymptotically optimal multistep predictors (or the order/method combinations with the minimal loss function values) for DGPs I–VIII are listed in Table 2. We generated 1000 replications for each of these DGPs and carried out predictor selection for each replication. The frequency of these combinations selected by Procedures A, B and C is shown in Table 3 for $n = 150, 300, 500, 1000$ and $2000$. The simulation results are summarized as follows:

(1) Two-step predictions. Procedures A, B and C can efficiently select the best order/method combination (listed in Table 2) regardless of whether the DGP is stationary or non-stationary. (Note that DGPs I and II are stationary, but DGPs III and IV are not.) In particular, the proportion of the best combination selected by Procedures B and C always exceeds 95 percent, except in DGPs II and IV with $n = 150$. Note that while the differences between the parameter values of DGPs I and II (or III and IV) are not sizable, different order/method combinations are required to attain the minimal loss function value (defined in (3.1) and (3.2)). Table 3 shows that these procedures are sensitive to small parameter changes and can efficiently switch to the "right track". However, we also notice that the finite-sample performance of Procedure A seems to be slightly worse than that of Procedures B and C.

(2) Three-step predictions. Note that DGPs V and VI are stationary AR(2) models with AR coefficients satisfying $0 < a_1 < 1$ and $a_1^2 + a_2 = 0$. Ing [10] recently showed that $(1, 2)$ is asymptotically more efficient than $(2, 1)$ in DGP V, whereas $(2, 1)$ is asymptotically more efficient than $(1, 2)$ in DGP VI. Procedures A, B and C perform quite well in this subtle case. More specifically, for $(a_1, a_2) = (0.9, -0.81)$, they can correctly choose $(1, 2)$ over 90 percent of the time for all sample sizes (except for Procedure A in the sample sizes of 150 and 300). On the other hand, when $(a_1, a_2) = (0.6, -0.36)$, Procedures B and C successfully select another combination, $(2,1)$, with rather high frequency for $n \geq 300$. While Procedure A performs slightly worse than the other two procedures, it can still choose $(2,1)$ with over 89 percent frequency as $n \geq 500$. Data generating processes VII and VIII are unit root processes. In DGP VII, the direct method only requires two regressors to perform three-step predictions and, according to Section 2, $(2, 2)$ can attain the minimal loss function value. On the other hand, $(3, 1)$ is the best combination for DGP VIII. Table 3 shows that the performance of Procedures A, B and C in DGPs VII and VIII are similar to those in DGPs V and VI.

**Table 2.** Order/method combination with the minimal loss function value

| DGP | $h = 2$ | | | | $h = 3$ | | | |
|---|---|---|---|---|---|---|---|---|
| | I | II | III | IV | V | VI | VII | VIII |
| Combination | (1,2) | (2,1) | (2,2) | (3,1) | (1,2) | (2,1) | (2,2) | (3,1) |

To explore the finite-sample performance of these procedures for larger lead times, we also conduct a small Monte Carlo study using the following two unit root AR models:

**DGP IX** $x_t = 0.2x_{t-10} + 0.8x_{t-11} + \varepsilon_t,$
**DGP X** $x_t = 1.5x_{t-1} - 0.5x_{t-2} + \varepsilon_t,$

where $\varepsilon_t$'s are independent and identically $\mathcal{N}(0, 25)$ distributed. Our goal is to select ten-step ($h = 10$) predictors for these two DGPs among a family of predictors, $\{(i, j), i = 1, \ldots, 20, j = 1, 2\}$. Note that DGP IX is an AR(11) model with $p_{10} = 2 \ll p_1 = 11$. Theorems 2.2 and 2.3 yield that the best combination for DGP IX is (2, 2). On the other hand, DGP X is an AR(2) model with $p_{10} = p_1 = 2$ and, in view of Example 1, (2, 1) is the best combination for DGP X. Our simulation results, based on 1000 replications for $n = 500$ and 1000, are reported in Table 4. Table 4 shows that when $h$ increases to 10, Procedures A, B and C still work well, except in DGP X with $n = 500$. In this latter case, while the proposed procedures can choose the best combination 70–80 percent of the time, we have found that the proportion of (1, 2) chosen by them is about 20 percent, indicating an underfitting problem. However, this difficulty is alleviated as $n$ increases to 1000, which coincides with the asymptotic results given in Theorem 4.1.

Finally, we note that the choice of $C_n$ in Procedure II does influence its finite-sample results. While we do not intend to suggest the best $C_n$ in finite-sample cases, the $C_n$'s used in this paper may serve as good "initial values" for pursuing better performance based on Procedure II.

**Table 3.** Frequency of choosing predictors with minimal loss function values in 1000 replications

| | | $h = 2$ | | | | $h = 3$ | | |
| | | Procedure | | | | Procedure | | |
| $n$ | Model(Unit Root) | A | B | C | Model(Unit Root) | A | B | C |
|---|---|---|---|---|---|---|---|---|
| 150 | I (No) | 853 | 963 | 987 | V (No) | 882 | 976 | 993 |
| 300 | | 890 | 984 | 997 | | 880 | 974 | 993 |
| 500 | | 901 | 990 | 999 | | 913 | 985 | 997 |
| 1000 | | 921 | 990 | 997 | | 918 | 994 | 999 |
| 2000 | | 948 | 992 | 1000 | | 951 | 991 | 1000 |
| 150 | II (No) | 817 | 887 | 869 | VI (No) | 698 | 711 | 689 |
| 300 | | 845 | 968 | 983 | | 827 | 936 | 915 |
| 500 | | 891 | 980 | 996 | | 898 | 989 | 992 |
| 1000 | | 913 | 985 | 995 | | 913 | 992 | 1000 |
| 2000 | | 923 | 990 | 999 | | 941 | 997 | 1000 |
| 150 | III (Yes) | 844 | 972 | 991 | VII (Yes) | 841 | 970 | 993 |
| 300 | | 893 | 989 | 997 | | 855 | 978 | 993 |
| 500 | | 916 | 992 | 998 | | 911 | 989 | 998 |
| 1000 | | 939 | 993 | 999 | | 917 | 995 | 999 |
| 2000 | | 950 | 997 | 1000 | | 939 | 997 | 1000 |
| 150 | IV (Yes) | 780 | 894 | 878 | VIII (Yes) | 633 | 722 | 705 |
| 300 | | 881 | 971 | 995 | | 835 | 901 | 903 |
| 500 | | 881 | 973 | 993 | | 888 | 973 | 975 |
| 1000 | | 906 | 980 | 994 | | 930 | 990 | 996 |
| 2000 | | 926 | 989 | 999 | | 944 | 994 | 1000 |

**Table 4.** Frequency of choosing predictors with minimal loss function values in 1000 replications

| | | | | | | | | |
|---|---|---|---|---|---|---|---|---|
| | | | $h = 10$ | | | | | |
| | | Procedure | | | | | Procedure | |
| $n$ | Model (Unit Root) | A | B | C | Model (Unit Root) | A | B | C |
| 500 | IX (Yes) | 967 | 1000 | 1000 | X (Yes) | 726 | 787 | 719 |
| 1000 | | 981 | 997 | 1000 | | 808 | 927 | 936 |

# Appendix A

Throughout this section, we only consider the case $k \geq 2$ (recall that $k$ denotes the order of the working AR model) because the results for the case $k = 1$ can be verified similarly. We start with some useful lemmas.

**Lemma A.1.** *Assume that $\{x_t\}$ satisfies model (1.1) with $\{\varepsilon_t\}$ obeying (2.1). Then, for any $q > 0$ and $k \geq p_1$,*

$$E\|\hat{R}_n^{-1}(k)\|^q = \mathrm{O}(1), \tag{A.1}$$

*where $\hat{R}_n(k)$ is defined after (2.2) and for a matrix $A$, $\|A\|^2 = \sup_{\|\mathbf{z}\|=1} \mathbf{z}'A'A\mathbf{z}$ with $\|\mathbf{z}\|$ denoting the Euclidean norm for vector $\mathbf{z}$.*

**Proof.** (A.1) can be verified by an argument similar to that used in the proof of Lemma A.1 in Ing *et al.* [13]. The details are omitted. □

**Lemma A.2.** *Assume that $\{x_t\}$ satisfies model (1.1) with $\{\varepsilon_t\}$ obeying (2.1) and for some $q_1 \geq 2$, $\sup_{-\infty < t < \infty} E|\varepsilon_t|^{2q_1} < \infty$. Then, for any $0 < q < q_1$ and $k \geq p + 1$,*

$$E\|\hat{R}_n^{-1}(k) - \hat{R}_n^{*-1}(k)\|^q = \mathrm{O}(n^{-q/2}), \tag{A.2}$$

*where*

$$\hat{R}_n^*(k) = \begin{pmatrix} \hat{\Gamma}_n(k-1) & \mathbf{0}'_{k-1} \\ \mathbf{0}'_{k-1} & \dfrac{1}{n^2}\displaystyle\sum_{j=k}^{n-1} N_j^2 \end{pmatrix},$$

*$\hat{\Gamma}_n(k-1) = (1/n)\sum_{j=k}^{n-1} \mathbf{s}_j(k-1)\mathbf{s}'_j(k-1)$ and $N_j = x_j - \sum_{l=1}^{k-1} \alpha_j x_{j-l}$.*

**Proof.** First note that Lemma A.1 ensures for any $q > 0$,

$$E\|\hat{R}_n^{*-1}(k)\|^q = \mathrm{O}(1). \tag{A.3}$$

We also have

$$\|\hat{R}_n^{-1}(k) - \hat{R}_n^{*-1}(k)\|^q \le \|\hat{R}_n^{-1}(k)\|^q \|\hat{R}_n^{*-1}(k)\|^q \|\hat{R}_n(k) - \hat{R}_n^*(k)\|^q$$

$$\le C_1 \|\hat{R}_n^{-1}(k)\|^q \|\hat{R}_n^{*-1}(k)\|^q \left\| n^{-3/2} \sum_{j=k}^{n-1} \mathbf{s}_j(k-1)N_j \right\|^q, \quad (A.4)$$

where $C_1$ is some positive constant. By analogy with Lemma A.3 in Ing *et al.* [13],

$$E \left\| n^{-3/2} \sum_{j=k}^{n-1} \mathbf{s}_j(k-1)N_j \right\|^{q_1} = \mathrm{O}(n^{-q_1/2}). \quad (A.5)$$

Consequently, (A.2) follows from (A.1), (A.3)–(A.5) and Hölder's inequality. $\square$

To prove Theorem 2.2, we also need the following two lemmas, the proofs of which are straightforward and hence omitted.

**Lemma A.3.** *Assume that $\{x_t\}$ satisfies model (1.1) with $\sup_{-\infty < t < \infty} E|\varepsilon_t|^q < \infty$, where $q \ge 2$. Then, for $k \ge p_1$,*

$$E \left\| n^{-1/2} D_n(k) \sum_{j=k}^{n-1} \mathbf{x}_j(k)\varepsilon_{j+1} \right\|^q = \mathrm{O}(1). \quad (A.6)$$

**Lemma A.4.** *Assume that $\{x_t\}$ satisfies model (1.1) with $\sup_{-\infty < t < \infty} E|\varepsilon_t|^r < \infty$ for some $r > 4$. Then, for $k \ge p_1$,*

$$\lim_{n \to \infty} E(F_{n,k}) = 0, \quad (A.7)$$

*where*

$$F_{n,k} = \frac{\mathbf{s}_n(k-1)M_h(k-1)\hat{\Gamma}_n^{-1}(k-1)\{\sum_{j=k}^{n-1} \mathbf{s}_j(k-1)\varepsilon_{j+1}\}N_n \sum_{j=k}^{n-1} N_j\varepsilon_{j+1}}{\sum_{j=k}^{n-1} N_j^2}. \quad (A.8)$$

**Proof of Theorem 2.2.** Some algebraic manipulations give

$$x_{n+h} - \hat{x}_{n+h}(k) = \eta_{n,h} - \mathbf{x}_n'(k)\hat{L}_{n,h}(k)(\hat{\mathbf{a}}_n(1,k) - \mathbf{a}(k)), \quad (A.9)$$

where $\hat{L}_{n,h}(k)$ is defined after (3.4). We also have

$$nE\{\mathbf{x}_n'(k)(\hat{L}_{n,h}(k) - L_h(k))(\hat{\mathbf{a}}_n(1,k) - \mathbf{a}(k))\}^2$$

$$= E\left\{ \mathbf{x}_n'(k)(\hat{L}_{n,h}(k) - L_h(k))D_n'(k)\hat{R}_n^{-1}(k)\frac{1}{\sqrt{n}}D_n(k) \sum_{j=k}^{n-1} \mathbf{x}_j(k)\varepsilon_{j+1} \right\}^2 \quad (A.10)$$

$$\equiv E\{G_n^2(k)\},$$

where $L_h(k) = \sum_{j=0}^{h-1} b_j A^{h-1-j}(k)$, with $A(k)$ defined in Section 1. Let $\hat{\alpha}_n(k-1) = (\hat{\alpha}(n,1), \ldots, \hat{\alpha}(n,k-1))' = \Psi(k)\hat{\mathbf{a}}_n(1,k)$, where $\Psi(k)$ is a $(k-1) \times k$ matrix, with the $(i,j)$th component equal to 0 if $j \leqq i$ and equal to $-1$ if $j > i$. Then, by observing $A(k)D'_n(k) = D'_n(k)\bar{A}(k)$ and $\hat{A}_n(k)\hat{D}'_n(k) = \hat{D}'_n(k)A^*_n(k)$, where $\hat{D}_n(k)$ is $D_n(k)$ with $\alpha_i$ replaced by $\hat{\alpha}(n,i)$ for $i = 1, \ldots, k-1$,

$$\bar{A}(k) = \begin{pmatrix} S_M(k-1) & \mathbf{0}_{k-1} \\ \mathbf{0}'_{k-1} & 1 \end{pmatrix},$$

and $A^*_n(k)$ is $\bar{A}(k)$ with $S_M(k-1)$ replaced by

$$\hat{S}_{M,n}(k-1) = \left( \hat{\alpha}_n(k-1) \,\middle|\, \begin{array}{c} I_{k-2} \\ \hline \mathbf{0}'_{k-2} \end{array} \right),$$

we have

$$L_h(k)D'_n(k) = D'_n(k)\bar{L}_h(k) \tag{A.11}$$

and

$$\hat{L}_{n,h}(k)\hat{D}'_n(k) = \hat{D}'_n(k)L^*_{n,h}(k), \tag{A.12}$$

where

$$\bar{L}_h(k) = \begin{pmatrix} M_h(k-1) & \mathbf{0}_{k-1} \\ \mathbf{0}'_{k-1} & \sum_{j=0}^{h-1} b_j \end{pmatrix}$$

and $L^*_{n,h}(k)$ is $\bar{L}_h(k)$ with $M_h(k-1)$ replaced by $\hat{M}_{n,h}(k-1) = \sum_{j=0}^{h-1} b_j \hat{S}^{h-1-j}_{M,n}(k-1)$. (A.11) and (A.12) yield

$$\begin{aligned}
(\hat{L}_{n,h}(k) &- L_h(k))D'_n(k) \\
&= \hat{L}_{n,h}(k)(D'_n(k) - \hat{D}'_n(k)) + (\hat{D}'_n(k) - D'_n(k))L^*_{n,h}(k) \\
&\quad + D'_n(k)(L^*_{n,h}(k) - \bar{L}_h(k)),
\end{aligned}$$

and hence

$$|G_n(k)| \leq G^*_n(k), \tag{A.13}$$

where $G^*_n(k) = (I) + (II)$, with

$$(I) = \|n^{-1/2}\mathbf{x}_n(k)\|\|\hat{\alpha}_n(k-1) - \alpha(k-1)\|(\|\hat{L}_{n,h}(k)\| + \|L^*_{n,h}(k)\|)G^*_{1,n}(k),$$

$$(II) = (\|\mathbf{s}_n(k-1)\| + |n^{-1/2}N_n|)\|L^*_{n,h}(k) - \bar{L}_h(k)\|G^*_{1,n}(k)$$

and $G^*_{1,n}(k) = \|\hat{R}^{-1}_n(k)\| \| n^{-1/2} D_n(k) \sum_{j=k}^{n-1} \mathbf{x}_j(k) \varepsilon_{j+1} \|$. By (A.10), (A.13), Lemmas A.1–A.3 and Hölder's inequality, it can be shown that

$$nE\{\mathbf{x}'_n(k)(\hat{L}_{n,h}(k) - L_h(k))(\hat{\mathbf{a}}_n(1,k) - \mathbf{a}(k))\}^2 \le E(G^*_n(k))^2 = \mathrm{O}(n^{-1}). \quad (A.14)$$

Similarly, we have

$$E\left\{\mathbf{x}'_n(k)L_h(k)D'_n(k)(\hat{R}^{-1}_n(k) - \hat{R}^{*-1}_n(k))\frac{1}{\sqrt{n}}D_n(k)\sum_{j=k}^{n-1}\mathbf{x}_j(k)\varepsilon_{j+1}\right\}^2 = \mathrm{O}(n^{-1}). \quad (A.15)$$

By (A.11) and some algebraic manipulations,

$$E\left\{\mathbf{x}'_n(k)L_h(k)D'_n(k)\hat{R}^{*-1}_n(k)\frac{1}{\sqrt{n}}D_n(k)\sum_{j=k}^{n-1}\mathbf{x}_j(k)\varepsilon_{j+1}\right\}^2$$
$$= E_{1,n}(k) + E_{2,n}(k) + E_{3,n}(k), \quad (A.16)$$

where

$$E_{1,n}(k) = E\left\{\mathbf{s}'_n(k-1)M_h(k-1)\hat{\Gamma}^{-1}_n(k-1)n^{-1/2}\sum_{j=k}^{n-1}\mathbf{s}_j(k-1)\varepsilon_{j+1}\right\}^2,$$

$$E_{2,n}(k) = \left(\sum_{j=0}^{h-1}b_j\right)^2 E\left\{n\frac{N^2_n(\sum_{j=k}^{n-1}N_j\varepsilon_{j+1})^2}{(\sum_{j=k}^{n-1}N^2_j)^2}\right\},$$

$$E_{3,n}(k) = 2\left(\sum_{j=0}^{h-1}b_j\right)E(F_{n,k}).$$

By an analogy with Theorem 1 of Ing [9],

$$\lim_{n\to\infty}E_{1,n}(k) = f_{1,h}(k-1). \quad (A.17)$$

In view of Ing [8], it is straightforward to show that

$$\lim_{n\to\infty}E_{2,n}(k) = 2\sigma^2\left(\sum_{j=0}^{h-1}b_j\right)^2. \quad (A.18)$$

Consequently, the desired result follows from (A.9), (A.14)–(A.18) and Lemma A.4. $\square$

**Proof of Theorem 2.3.** By analogies with lemmas A.1–A.4, for $k \ge p_h$,

$$E\|\bar{R}^{-1}_{n,h}(k)\|^q = \mathrm{O}(1), \quad (A.19)$$

$$E\|\bar{R}^{-1}_{n,h}(k) - \bar{R}^{*-1}_{n,h}(k)\|^4 = \mathrm{O}(n^{-2}), \quad (A.20)$$

$$E\left\|n^{-1/2}\bar{D}_n\sum_{j=k}^{n-1}\mathbf{x}_j(k)\eta_{j+h}\right\|^8 = \mathrm{O}(1) \tag{A.21}$$

and

$$\lim_{n\to\infty}E(\bar{F}_{n,k}) = 0, \tag{A.22}$$

where $q > 0$,

$$\bar{R}_{n,h}^*(k) = \begin{pmatrix} \dfrac{1}{n}\sum_{j=k}^{n-h}\mathbf{s}_j(k-1)\mathbf{s}_j'(k-1) & \mathbf{0}_{k-1} \\[2ex] \mathbf{0}_{k-1}' & \dfrac{1}{n^2}\sum_{j=k}^{n-h}x_j^2 \end{pmatrix}$$

and

$$\bar{F}_{n,k} = \frac{\mathbf{s}_n(k-1)\{(1/n)\sum_{j=k}^{n-h}\mathbf{s}_j(k-1)\mathbf{s}_j'(k-1)\}^{-1}\{\sum_{j=k}^{n-h}\mathbf{s}_j(k-1)\eta_{j,h}\}x_n\sum_{j=k}^{n-h}x_j\eta_{j,h}}{\sum_{j=k}^{n-h}x_j^2}. \tag{A.23}$$

In addition, according to (1.9) and (3.5) of Ing and Sin [12], it can be shown that

$$\lim_{n\to\infty}E\left\{n\frac{x_n^2(\sum_{j=k}^{n-h}x_j\eta_{j,h})^2}{\sum_{j=k}^{n-h}x_j^2}\right\} = 2\sigma^2\left(\sum_{j=0}^{h-1}b_j\right)^2. \tag{A.24}$$

As a result, Theorem 2.3 follows from (A.19)–(A.22), (A.24) and arguments similar to those used in the proofs of Theorem 2 in Ing [9] and Theorem 2.2 above.  □

# Appendix B

Lemma B.1 below provides (almost sure) asymptotic bounds for $\|\hat{\Gamma}_n(k-1) - \Gamma(k-1)\|$, $\|\hat{R}_n(k) - \hat{R}_n^*(k)\|$ and $\|\hat{R}_n^{-1}(k)\|$ under a minimal moment condition, $\sup_{-\infty<t<\infty}E|\varepsilon_t|^\alpha$ for some $\alpha > 2$. As will be seen later, these bounds play subtle roles in our asymptotic analysis.

**Lemma B.1.** *Assume that the assumptions of Theorem 3.1 hold. Then,*

(i) *for $k \geq 2$, $k \geq p_1$, and some $\iota > 0$,*

$$\|\hat{\Gamma}_n(k-1) - \Gamma(k-1)\| = \mathrm{o}(n^{-\iota}) \qquad a.s.; \tag{B.1}$$

(ii) *for $k \geq p_1$ and some $\eta > 0$,*

$$\|\hat{R}_n(k) - \hat{R}_n^*(k)\| = \mathrm{o}(n^{-\eta}) \qquad a.s.; \tag{B.2}$$

(iii) *for $k \geq p_1$,*

$$\|\hat{R}_n^{-1}(k)\| = \mathrm{O}(\log \log n) \qquad a.s. \tag{B.3}$$

**Proof.** First note that

$$\|\hat{\Gamma}_n(k-1) - \Gamma(k-1)\| \leq \sum_{l=0}^{k-2} \sum_{m=0}^{k-2} \left| n^{-1} \sum_{j=k}^{n-1} s_{j-l} s_{j-m} - \gamma_{l,m} \right|,$$

where $\gamma_{l,m}$ is the $(l,m)$th component of $\Gamma(k-1)$. Therefore, (B.1) is ensured by showing that for any $1 \leq l \leq k-1$ and $1 \leq m \leq k-1$,

$$\left| \frac{1}{n} \sum_{j=k}^{n-1} s_{j-l} s_{j-m} - \gamma_{l,m} \right| = \mathrm{o}(n^{-\iota}) \qquad a.s. \tag{B.4}$$

In the following, we only prove the case of $l = m = 0$ since the proofs of other cases can be similarly obtained. For $l = m = 0$, the left-hand side of (B.4) can be rewritten as

$$\left| \frac{1}{n} \sum_{j=k}^{n-1} (s_j^2 - \gamma_{0,0}^{(j)}) + \frac{1}{n} \sum_{j=k}^{n-1} (\gamma_{0,0}^{(j)} - \gamma_{0,0}) + \frac{k\gamma_{0,0}}{n} \right|, \tag{B.5}$$

where $\gamma_{0,0}^{(j)} = \sigma^2 \sum_{r=0}^{j-1} c_r^2$ with $c_j$'s defined in Section 1. By observing $\gamma_{0,0} = \sigma^2 \sum_{r=0}^{\infty} c_r^2$ and $|c_r| \leq C_1 \mathrm{e}^{-\beta_1 r}$ for all $r$ and some $C_1, \beta_1 > 0$, we have $(1/n) \sum_{j=k}^{n-1} (\gamma_{0,0}^{(j)} - \gamma_{0,0}) = \mathrm{O}(1/n)$ and $k\gamma_{0,0}/n = \mathrm{O}(1/n)$. In addition, straightforward calculations yield that

$$s_j^2 - \gamma_{0,0}^{(j)} = \sum_{l=1}^{j} c_{j-l}^2 (\varepsilon_l^2 - \sigma^2) + 2 \sum_{l_2=2}^{j} \sum_{l_1=1}^{l_2-1} c_{j-l_1} c_{j-l_2} \varepsilon_{l_1} \varepsilon_{l_2}. \tag{B.6}$$

In view of (B.6), one obtains, through changing the order of summations, that

$$\begin{aligned}
\sum_{j=n_1}^{n_2} \frac{s_j^2 - \gamma_{0,0}^{(j)}}{j^{\theta}} &= \sum_{l=1}^{n_1} \left( \sum_{j=n_1}^{n_2} \frac{c_{j-l}^2}{j^{\theta}} \right) \eta_l + \sum_{l=n_1+1}^{n_2} \left( \sum_{j=l}^{n_2} \frac{c_{j-l}^2}{j^{\theta}} \right) \eta_l \\
&\quad + 2 \sum_{l_2=2}^{n_1} \left\{ \sum_{l_1=1}^{l_2-1} \left( \sum_{j=n_1}^{n_2} \frac{c_{j-l_1} c_{j-l_2}}{j^{\theta}} \right) \varepsilon_{l_1} \right\} \varepsilon_{l_2} \\
&\quad + 2 \sum_{l_2=n_1+1}^{n_2} \left\{ \sum_{l_1=1}^{l_2-1} \left( \sum_{j=l_2}^{n_2} \frac{c_{j-l_1} c_{j-l_2}}{j^{\theta}} \right) \varepsilon_{l_1} \right\} \varepsilon_{l_2} \equiv (I) + (II) + (III) + (IV),
\end{aligned}$$

where $\eta_t = \varepsilon_t^2 - \sigma^2$, $\theta < 1$ and $\theta\alpha/2 > 1$. If we can show that for any $1 \le n_1 \le n_2 < \infty$,

$$E|(G)|^{\alpha/2} \le C\left(\sum_{j=n_1}^{n_2} \frac{1}{j^{\xi_1}}\right)^{\xi_2} \tag{B.7}$$

(where $G = I, II, III$ and $IV$) and $C > 0, \xi_1 > 1, \xi_2 > 1$ are some positive constant independent of $n_1$ and $n_2$ (but they can vary with $G$), then by Móricz [18] for all sufficiently large $n_1$,

$$E \max_{n_1 \le l \le n_2} \left|\sum_{j=n_1}^{l} \frac{s_j^2 - \gamma_{0,0}^{(j)}}{j^\theta}\right|^{\alpha/2} \le C^*\left(\sum_{j=n_1}^{n_2} \frac{1}{j^{\xi_1^*}}\right)^{\xi_2^*}, \tag{B.8}$$

where $C^* > 0$, $\xi_1^* > 1$ and $\xi_2^* > 1$ are some positive constants independent of $n_1$ and $n_2$. (B.8) and Kronecker's lemma yield

$$\frac{1}{n^\theta}\sum_{j=1}^{n}(s_j^2 - \gamma_{0,0}^{(j)}) = o(1) \qquad \text{a.s.} \tag{B.9}$$

As a result, (B.1) holds with $\iota = 1 - \theta$.

Without loss of generality, assume $2 < \alpha < 4$. Then,

$$\begin{aligned}
E|(I)|^{\alpha/2} &\le C_2 E\left\{\sum_{l=1}^{n_1}\left(\sum_{j=n_1}^{n_2} \frac{c_{j-l}^2}{j^\theta}\right)^2 \eta_l^2\right\}^{\alpha/4} \\
&\le C_2 \sum_{j_1=n_1}^{n_2}\sum_{j_2=n_1}^{n_2} \frac{1}{j_1^{\theta\alpha/4}j_2^{\theta\alpha/4}}\sum_{l=1}^{n_1}|c_{j_1-l}c_{j_2-l}|^{\alpha/2}E|\eta_l|^{\alpha/2} \\
&\le C_3\left(\sum_{j=n_1}^{n_2}\frac{1}{j^{\theta\alpha/2}} + \sum_{j_1=n_1}^{n_2-1}\frac{1}{j_1^{\theta\alpha/4}}\sum_{j_2=j_1+1}^{n_2}\frac{1}{j_2^{\theta\alpha/4}}(j_2-j_1)^{-s}\right) \\
&\le C_4\left(\sum_{j=n_1}^{n_2}\frac{1}{j^{\theta\alpha/2}}\right) \\
&\le C_4\left(\sum_{j=n_1}^{n_2}\frac{1}{j^{\xi_1}}\right)^{\xi_2},
\end{aligned} \tag{B.10}$$

where $C_i > 0, i = 2, \ldots, 4$, and $s > 1$ are some positive constants independent of $n_1$ and $n_2$, $1 < \xi_1 < \theta\alpha/2$, $\xi_2 = \theta\alpha/2\xi_1$, the first inequality follows from Burkholder's inequality, the second one follows from the fact that $\alpha/4 < 1$ and changing the order of summations, the third one is ensured by $\sup_t E|\varepsilon_t|^\alpha < \infty$ and $c_j \le C_1 e^{-\beta_1 j}$, which implies for all $n_1 \le j_1 \ne j_2 \le n_2$, $\sum_{l=1}^{n_1}|c_{j_1-l}c_{j_2-l}|^{\alpha/2} \le C_5|j_1-j_2|^{-s}$, for some $C_5 > 0$. As a result, (B.7) holds with $G = I$. The proof of (B.7) for the case of $G = II$ is similar. The details are omitted. To show (B.7) for the case of $G = III$, note that

$$
\begin{aligned}
E\Bigg| &\sum_{l_2=2}^{n_1} \Bigg\{ \sum_{l_1=1}^{l_2-1} \Bigg( \sum_{j=n_1}^{n_2} \frac{c_{j-l_1} c_{j-l_2}}{j^\theta} \Bigg) \varepsilon_{l_1} \Bigg\} \varepsilon_{l_2} \Bigg|^{\alpha/2} \\
&\leq \Bigg\{ E \Bigg( \sum_{l_2=2}^{n_1} \Bigg\{ \sum_{l_1=1}^{l_2-1} \Bigg( \sum_{j=n_1}^{n_2} \frac{c_{j-l_1} c_{j-l_2}}{j^\theta} \Bigg) \varepsilon_{l_1} \Bigg\} \varepsilon_{l_2} \Bigg)^2 \Bigg\}^{\alpha/4} \\
&= |\sigma|^\alpha \Bigg( \sum_{l_2=2}^{n_1} \sum_{l_1=1}^{l_2-1} \Bigg( \sum_{j=n_1}^{n_2} \frac{c_{j-l_1} c_{j-l_2}}{j^\theta} \Bigg)^2 \Bigg)^{\alpha/4}.
\end{aligned}
\tag{B.11}
$$

By arguments similar to those used to verify the second to fifth inequalities in (B.10), the desired result follows. Similarly, it can be shown that (B.7) holds for the case of $G = IV$.

To show (B.2), first observe that

$$
\|\hat{R}_n(k) - \hat{R}_n^*(k)\| \leq \sqrt{2} \sum_{l=0}^{k-2} \Bigg| \frac{1}{n^{3/2}} \sum_{j=k}^{n-1} s_{j-l} N_j \Bigg|.
$$

Therefore, it suffices to show that for $l = 0, \ldots, k-2$ and some $\eta > 0$,

$$
\frac{1}{n^{3/2}} \sum_{j=k}^{n-1} s_{j-l} N_j = \mathrm{o}(n^{-\eta}) \qquad \text{a.s.}
\tag{B.12}
$$

We only verify (B.12) for the case $l = 0$ since the proof of the case $l > 0$ can be similarly obtained. Let $\max\{1, (1/2) + (2/\alpha)\} < \theta_1 < 3/2$. Some algebraic manipulations yield

$$
\begin{aligned}
\sum_{j=n_1}^{n_2} \frac{s_j N_j}{j^{\theta_1}} = {}& \sigma^2 \sum_{j=n_1}^{n_2} \frac{1}{j^{\theta_1}} \sum_{m=1}^{j} c_{j-m} + \sum_{m=1}^{n_1} \sum_{j=n_1}^{n_2} \frac{c_{j-m}}{j^{\theta_1}} \eta_m + \sum_{m=n_1+1}^{n_2} \sum_{j=m}^{n_2} \frac{c_{j-m}}{j^{\theta_1}} \eta_m \\
&+ \sum_{m=2}^{n_1} \sum_{l=1}^{m-1} \sum_{j=n_1}^{n_2} \frac{c_{j-l}}{j^{\theta_1}} \varepsilon_l \varepsilon_m + \sum_{m=n_1+1}^{n_2} \sum_{l=1}^{m-1} \sum_{j=m}^{n_2} \frac{c_{j-l}}{j^{\theta_1}} \varepsilon_l \varepsilon_m \\
&+ \sum_{l=2}^{n_1} \sum_{m=1}^{l-1} \sum_{j=n_1}^{n_2} \frac{c_{j-l}}{j^{\theta_1}} \varepsilon_m \varepsilon_l + \sum_{l=n_1+1}^{n_2} \sum_{m=1}^{l-1} \sum_{j=l}^{n_2} \frac{c_{j-l}}{j^{\theta_1}} \varepsilon_m \varepsilon_l \\
= {}& (\ddot{I}) + (\ddot{II}) + (\ddot{III}) + (\ddot{IV}) + (\ddot{V}) + (\ddot{VI}) + (\ddot{VII}).
\end{aligned}
\tag{B.13}
$$

It is clear that

$$
|(\ddot{I})|^{\alpha/2} \leq C_6 \Bigg( \sum_{j=n_1}^{n_2} \frac{1}{j^{\xi_1}} \Bigg)^{\xi_2},
\tag{B.14}
$$

where $\xi_1 = \theta_1$ and $\xi_2 = \alpha/2$. By an argument similar to that used in (B.10),

$$E|(W)|^{\alpha/2} \le C_7 \left( \sum_{j=n_1}^{n_2} \frac{1}{j^{\xi_1}} \right)^{\xi_2}, \tag{B.15}$$

where $W = \ddot{I}\ddot{I}, I\ddot{I}\ddot{I}$, $1 < \xi_1 < \theta_1\alpha/2$ and $\xi_2 = \theta_1\alpha/2\xi_1$. An argument similar to that used in (B.11) yields

$$E|(W)|^{\alpha/2} \le C_8 \left( \sum_{j=n_1}^{n_2} \frac{1}{j^{\xi_1}} \right)^{\xi_2}, \tag{B.16}$$

where $W = I\ddot{V}, \ddot{V}, \ddot{V}I, \ddot{V}II$, $1 < \xi_1 < (2\theta_1-1)\alpha/4$ and $\xi_2 = (2\theta_1-1)\alpha/4\xi_1$. Consequently, (B.12) (with $\eta = (3/2) - \theta_1$) follows from (B.13)–(B.16), Móricz [18] and Kronecker's lemma. To show (B.3), observe that $\|\hat{R}_n^{-1}(k)\| \le \|\hat{R}_n^{-1}(k)\| \|\hat{R}_n(k) - \hat{R}_n^*(k)\| \|\hat{R}_n^{*-1}(k)\| + \|\hat{R}_n^{*-1}(k)\|$. By (3.23) of Lai and Wei [14] and (3.2) of Lai and Wei [15],

$$\|\hat{R}_n^{*-1}(k)\| = \mathrm{O}(\log\log n) \qquad \text{a.s.}$$

This and (B.2) yield (B.3). □

To prove Theorem 3.1, the following auxiliary lemma is required. Its proof can be found in Appendix B of Ing *et al.* [11].

**Lemma B.2.** *Assume that the assumptions of Theorem 3.1 hold. Then, for $k \ge \max\{2, p_1\}$,*

$$\sum_{i=m_h}^{n-h} F_{i,k} = \mathrm{o}(n) \qquad a.s., \tag{B.17}$$

*where $F_{i,k}$ is defined in Lemma A.4.*

We also need a few elementary facts.

**Lemma B.3.** *Let $\{z_n\}$ be a sequence of real numbers.*

(i) *If $z_n \ge 0$, $n^{-1}\sum_{j=1}^n z_j = \mathrm{O}(1)$ and, for some $\xi > 1$, $\liminf_{n\to\infty} \nu_n/n^\xi > 0$, then*

$$\sum_{j=1}^n \frac{z_j}{\nu_j} = \mathrm{O}(1).$$

(ii) *If $n^{-1}\sum_{j=1}^n z_j = \mathrm{o}(1)$, then*

$$\sum_{j=1}^n \frac{z_j}{j} = \mathrm{o}(\log n).$$

**Proof of Theorem 3.1.** We only prove the case $k \geq 2$ since the proof of the case $k = 1$ is similar. By Chow [1] and an analogy with (3.8) of Ing [10],

$$\text{APE}P_{n,h}(k) - \sum_{i=m_h}^{n-h} (\eta_{i,h})^2 = \sum_{i=m_h}^{n-h} \{\mathbf{x}_i'(k)\hat{L}_{i,h}(k)(\hat{\mathbf{a}}_i(1,k) - \mathbf{a}(k))\}^2 (1 + \text{o}(1))$$
$$+ \text{O}(1) \quad \text{a.s.} \tag{B.18}$$

Straightforward calculations give

$$\sum_{i=m_h}^{n-h} \{\mathbf{x}_i'(k)(\hat{L}_{i,h}(k) - L_h(k))(\hat{\mathbf{a}}_i(1,k) - \mathbf{a}(k))\}^2$$
$$= \sum_{i=m_h}^{n-h} \frac{1}{i} \left\{ \mathbf{x}_i'(k)(\hat{L}_{i,h}(k) - L_h(k))D_i'(k)\hat{R}_i^{-1}(k)\frac{1}{\sqrt{i}}D_i(k)\sum_{j=k}^{i-1}\mathbf{x}_j(k)\varepsilon_{j+1} \right\}^2. \tag{B.19}$$

By Lai and Wei [14] and (3.1) and (3.2) of Lai and Wei [15], we have

$$\|\hat{\alpha}_n(k-1) - \alpha(k-1)\| = \text{O}\left(\left(\frac{\log n}{n}\right)^{1/2}\right) \quad \text{a.s.,} \tag{B.20}$$

$$\|L_{n,h}^*(k) - \bar{L}_h(k)\| = \text{O}\left(\left(\frac{\log n}{n}\right)^{1/2}\right) \quad \text{a.s.,} \tag{B.21}$$

$$\|\hat{L}_{n,h}(k)\| = \text{O}(1) \quad \text{a.s.,} \tag{B.22}$$

$$\|\mathbf{x}_n(k)/\sqrt{n}\| = \text{O}((\log\log n)^{1/2}) \quad \text{a.s.} \tag{B.23}$$

In addition, by Lemma 1 of Wei [21], the law of the iterated logarithm, and (3.3) of Lai and Wei [15],

$$\left\| \frac{1}{\sqrt{n}}D_n(k)\sum_{j=k}^{n-1}\mathbf{x}_j(k)\varepsilon_{j+1} \right\| = \text{o}((\log n)^{\delta}(\log\log n)^{1/2}) \quad \text{a.s.,} \tag{B.24}$$

where $\delta > 1/\alpha$. As a result, by (A.11), (A.12), (B.3), (B.19)–(B.24) and the fact that $N_n/\sqrt{n} = \text{O}((\log\log n)^{1/2})$ a.s., one obtains

$$\sum_{i=m_h}^{n-h} \{\mathbf{x}_i'(k)(\hat{L}_{i,h}(k) - L_h(k))(\hat{\mathbf{a}}_i(1,k) - \mathbf{a}(k))\}^2 = \text{O}(1) \quad \text{a.s.} \tag{B.25}$$

Armed with (B.2), (B.3) and the fact that $\|\hat{R}_n^{*-1}(k)\| = \text{O}(\log\log n)$ a.s. (which is given after (B.16)), it can be shown that

$$\|\hat{R}_n^{-1}(k) - \hat{R}_n^{*-1}(k)\| = \text{o}\left(\frac{(\log\log n)^2}{n^{\eta}}\right) \quad \text{a.s.,} \tag{B.26}$$

where $\eta > 0$ is some positive constant. Since (A.11) yields for some $C_1 > 0$, $\|D_i(k)L'_h(k) \times \mathbf{x}_i(k)\| = \|\tilde{L}'_h(k)D_i(k)\mathbf{x}_i(k)\| \leq C_1(\|\mathbf{s}_i(k-1)\| + |N_i/\sqrt{i}|)$, we obtain

$$\sum_{i=m_h}^{n-h} \frac{1}{i} \left\{ \mathbf{x}'_i(k)L_h(k)D'_i(k)(\hat{R}_i^{-1}(k) - \hat{R}_i^{*-1}(k)) \frac{1}{\sqrt{i}} D_i(k) \sum_{j=k}^{i-1} \mathbf{x}_j(k)\varepsilon_{j+1} \right\}^2$$

$$\leq C_2 \sum_{i=m_h}^{n-h} \frac{1}{i} \left\{ \left( \|\mathbf{s}_i(k-1)\| + \left| \frac{N_i}{\sqrt{i}} \right| \right) \|\hat{R}_i^{-1}(k) - \hat{R}_i^{*-1}(k)\| \right. \tag{B.27}$$

$$\left. \times \left\| \frac{1}{\sqrt{i}} D_i(k) \sum_{j=k}^{i-1} \mathbf{x}_j(k)\varepsilon_{j+1} \right\| \right\}^2 = \mathrm{O}(1) \qquad \text{a.s.,}$$

where $C_2 > 0$ is some positive constant independent of $n$ and the equality follows from (B.24), (B.26), $N_n/\sqrt{n} = \mathrm{O}((\log\log n)^{1/2})$ a.s., $(1/n)\sum_{j=k}^{n-1} \|\mathbf{s}_j(k-1)\| = \mathrm{O}(1)$ a.s. and (i) of Lemma B.3.

By (A.11) and some algebraic manipulations,

$$\sum_{i=m_h}^{n-h} \frac{1}{i} \left\{ \mathbf{x}'_i(k)L_h(k)D'_i(k)\hat{R}_i^{*-1}(k) \frac{1}{\sqrt{i}} D_i(k) \sum_{j=k}^{i-1} \mathbf{x}_j(k)\varepsilon_{j+1} \right\}^2 = (I) + (II) + (III),$$

where

$$(I) = \sum_{i=m_h}^{n-h} \left\{ \mathbf{s}'_i(k-1)M_h(k-1)\hat{\Gamma}_i^{-1}(k) \frac{1}{i} \sum_{j=k}^{i-1} \mathbf{s}_j(k-1)\varepsilon_{j+1} \right\}^2,$$

$$(II) = \left( \sum_{j=0}^{h-1} b_j \right)^2 \sum_{i=m_h}^{n-h} \frac{N_i^2(\sum_{j=k}^{i-1} N_j\varepsilon_{j+1})^2}{(\sum_{j=k}^{i-1} N_j^2)^2},$$

$$(III) = 2\left( \sum_{j=0}^{h-1} b_j \right) \sum_{i=m_h}^{n-h} \frac{F_{i,k}}{i}.$$

According to (B.21) and analogies with (A.1) and Theorem 3.1 of Ing [10],

$$(I) = \sum_{i=m_h}^{n-h} \left\{ \mathbf{s}'_i(k-1)\hat{M}_{i,h}(k-1)\hat{\Gamma}_i^{-1}(k) \frac{1}{i} \sum_{j=k}^{i-1} \mathbf{s}_j(k-1)\varepsilon_{j+1} \right\}^2 + \mathrm{o}(\log n) \qquad \text{a.s.}$$

$$= f_{1,h}(k-1)\log n + \mathrm{o}(\log n) \qquad \text{a.s.}$$

By Theorem 4 of Wei [21],

$$(II) = 2\left( \sum_{j=0}^{h-1} b_j \right)^2 \sigma^2 \log n + \mathrm{o}(\log n) \qquad \text{a.s.}$$

In view of Lemma B.2 and (ii) of Lemma B.3, one obtains

$$(III) = \text{o}(\log n) \qquad \text{a.s.}$$

As a result,

$$\sum_{i=m_h}^{n-h} \frac{1}{i} \left\{ \mathbf{x}_i'(k) L_h(k) D_i'(k) \hat{R}_i^{*^{-1}}(k) \frac{1}{\sqrt{i}} D_i(k) \sum_{j=k}^{i-1} \mathbf{x}_j(k) \varepsilon_{j+1} \right\}^2$$

$$= \left\{ 2 \left( \sum_{j=0}^{h-1} b_j \right)^2 \sigma^2 + f_{1,h}(k-1) \right\} \log n + \text{o}(\log n) \qquad \text{a.s.}$$ (B.28)

Consequently, (3.6) follows from (B.18), (B.25), (B.27), (B.28) and the Cauchy–Schwarz inequality. □

To analyze APED$_{n,h}(k)$, Lemma B.4 is required.

**Lemma B.4.** *Let the assumptions of Theorem 3.1 hold. Then,*

$$\sum_{i=m_h}^{n-h} \frac{x_i^2 (\sum_{j=k}^{i-h} x_j \eta_{j,h})^2}{(\sum_{j=k}^{i-h} x_j^2)^2} = 2 \left( \sum_{j=0}^{h-1} b_j \right)^2 \sigma^2 \log n + \text{o}(\log n) \qquad a.s.$$ (B.29)

**Proof.** Following arguments similar to those used in the proofs of Lemma 2 and Theorem 1 of Ing and Sin [12], one obtains

$$\liminf_{n \to \infty} \frac{\log \log n}{n^2} \sum_{j=1}^{n-1} x_j^2 > 0 \qquad \text{a.s.}$$ (B.30)

and

$$x_n = \text{O}((n \log \log n)^{1/2}) \qquad \text{a.s.}$$ (B.31)

By the Borel–Cantelli lemma,

$$\varepsilon_n = \text{o}(n^{1/2}) \qquad \text{a.s.}$$ (B.32)

In addition, it is not difficult to show that for $\theta > 1/2$ and $l \geq 1$,

$$\frac{1}{n^\theta} \sum_{j=1}^{n-l} \varepsilon_j \varepsilon_{j+l} = \text{o}(1) \qquad \text{a.s.}$$ (B.33)

(B.30)–(B.33) together imply

$$\sum_{i=m_h}^{n-h} \frac{x_i^2 (\sum_{j=k}^{i-h} x_j \eta_{j,h})^2}{(\sum_{j=k}^{i-h} x_j^2)^2} = \left( \sum_{j=0}^{h-1} b_j \right)^2 \sum_{i=m_h}^{n-h} \frac{x_i^2 (\sum_{j=k}^{i-1} x_j \varepsilon_{j+1})^2}{(\sum_{j=k}^{i-1} x_j^2)^2} + \text{O}(1) \qquad \text{a.s.}$$ (B.34)

Consequently, (B.29) follows from (B.34) and the fact that

$$\sum_{i=m_h}^{n-h} \frac{x_i^2 (\sum_{j=k}^{i-1} x_j \varepsilon_{j+1})^2}{(\sum_{j=k}^{i-1} x_j^2)^2} = 2\sigma^2 \log n + \mathrm{o}(\log n) \qquad \text{a.s.,}$$

which is guaranteed by (2.15) of Ing and Sin [12].                                    □

**Proof of Theorem 3.2.** We only prove the case of $k \geq 2$ since the proof of the case of $k = 1$ is similar. By the same reasoning as in (B.18), we have

$$\text{APE}D_{n,h}(k) - \sum_{i=m_h}^{n-h} \eta_{i,h}^2 = (1 + \mathrm{o}(1)) \sum_{i=m_h}^{n-h} \left\{ \mathbf{x}_i'(k) V_{i-h}(k) \left( \sum_{j=k}^{i-h} \mathbf{x}_j(k) \eta_{j,h} \right) \right\}^2$$
$$+ \mathrm{O}(1) \qquad \text{a.s.} \tag{B.35}$$

Observe that

$$\sum_{i=m_h}^{n-h} \left\{ \mathbf{x}_i'(k) V_{i-h}(k) \left( \sum_{j=k}^{i-h} \mathbf{x}_j(k) \eta_{j,h} \right) \right\}^2$$
$$= \sum_{i=m_h}^{n-h} \frac{1}{i} \left\{ \mathbf{x}_i'(k) \bar{D}_i'(k) \bar{R}_{i,h}^{-1}(k) \left( \frac{1}{\sqrt{i}} \bar{D}_i(k) \sum_{j=k}^{i-h} \mathbf{x}_j(k) \eta_{j,h} \right) \right\}^2. \tag{B.36}$$

According to (B.30), (B.31) and arguments similar to those used to obtain (B.24) and (B.26),

$$\left\| \frac{1}{\sqrt{n}} \bar{D}_n(k) \sum_{j=k}^{n-h} \mathbf{x}_j(k) \eta_{j,h} \right\| = \mathrm{o}((\log n)^\delta (\log \log n)^{1/2}) \qquad \text{a.s.,}$$

and

$$\| \bar{R}_{n,h}^{-1}(k) - \bar{R}_{n,h}^{*-1}(k) \| = \mathrm{O}(n^{-\eta}(\log \log n)^2) \qquad \text{a.s.,}$$

where $\delta > 1/\alpha$ and $\eta > 0$. These facts and reasoning similar to that used in (B.27) yield

$$\sum_{i=m_h}^{n-h} \frac{1}{i} \left\{ \mathbf{x}_i'(k) \bar{D}_i'(k) (\bar{R}_{i,h}^{-1}(k) - \bar{R}_{i,h}^{*-1}(k)) \left( \frac{1}{\sqrt{i}} \bar{D}_i(k) \sum_{j=k}^{i-h} \mathbf{x}_j(k) \eta_{j,h} \right) \right\}^2$$
$$= \mathrm{O}(1) \qquad \text{a.s.} \tag{B.37}$$

Now,

$$\sum_{i=m_h}^{n-h} \frac{1}{i} \left\{ \mathbf{x}_i'(k) \bar{D}_i'(k) \bar{R}_{i,h}^{*-1}(k) \left( \frac{1}{\sqrt{i}} \bar{D}_i(k) \sum_{j=k}^{i-h} \mathbf{x}_j(k) \eta_{j,h} \right) \right\}^2$$
$$= (I) + (II) + (III), \tag{B.38}$$

where

$$(I) = \sum_{i=m_h}^{n-h} \left\{ \mathbf{s}_i'(k-1) \left[ \sum_{j=k}^{i-h} \mathbf{s}_j(k-1)\mathbf{s}_j'(k-1) \right]^{-1} \sum_{j=k}^{i-h} \mathbf{s}_j(k-1)\eta_{j,h} \right\}^2,$$

$$(II) = \sum_{i=m_h}^{n-h} \frac{x_i^2 (\sum_{j=k}^{i-h} x_j \eta_{j,h})^2}{(\sum_{j=k}^{i-h} x_j^2)^2}$$

and

$$(III) = \sum_{j=m_h}^{n-h} \frac{\bar{F}_i(k)}{i}$$

with $\bar{F}_i(k)$ defined in (A.23). By analogy with Theorem 3.2 of Ing [10],

$$(I) = f_{2,h}(k-1) \log n + o(\log n) \qquad \text{a.s.} \tag{B.39}$$

According to Lemma B.4,

$$(II) = 2\sigma^2 \left( \sum_{j=0}^{h-1} b_j \right)^2 \log n + o(\log n) \qquad \text{a.s.} \tag{B.40}$$

By reasoning similar to that used in the proof of Lemma B.2 (see Appendix B of Ing *et al.* [11]),

$$\sum_{j=m_h}^{n-h} \bar{F}_i(k) = o(n) \qquad \text{a.s.,}$$

and hence

$$(III) = o(\log n) \qquad \text{a.s.} \tag{B.41}$$

Consequently, (3.8) follows from (B.35)–(B.41). $\qquad \square$

To prove Theorem 3.3 we need a technical lemma, the proof of which can also be found in Appendix B of Ing *et al.* [11].

**Lemma B.5.** *Let the assumptions of Theorem 3.1 hold. Then, for $1 \le k < p_h$ and $h \ge 1$,*

$$\sum_{i=m_h}^{n-h} \left[ \mathbf{x}_i'(k)(\check{\mathbf{a}}_i(h,k) - \check{\mathbf{a}}(h,k)) \right]^2 = o(n) \qquad \text{a.s.,} \tag{B.42}$$

*where* $\tilde{\mathbf{a}}(h,k) = 1$, *if* $1 = k < p_h$, *and*

$$\tilde{\mathbf{a}}(h,k) = U_{k \times (k-1)} \left\{ \sum_{j=0}^{h-1} \alpha_{h-j}(k-1) \right\} + (1,0,\ldots,0)', \tag{B.43}$$

*if* $1 < k < p_h$, *where* $U_{k \times (k-1)} = (u_{ij})$ *is a* $k \times (k-1)$ *matrix, with* $u_{ij} = 1$, *if* $i = j$, $u_{ij} = -1$, *if* $i - j = 1$ *and* $u_{ij} = 0$, *otherwise, and* $\alpha_l(k-1) = \lim_{t \to \infty} \alpha_l^{(t)}(k-1)$, *with*

$$\alpha_l^{(t)}(k-1) = \arg \min_{(f_1,\ldots,f_{k-1})' \in R^{k-1}} E(s_{t+l} - f_1 s_t - \cdots - f_{k-1} s_{t-k+2})^2.$$

**Proof of Theorem 3.3.** By an analogy with (B.35),

$$\begin{aligned}
\mathrm{APED}&_{n,h}(k) \\
&= \sum_{i=m_h}^{n-h} \{\eta_{i,h} + \mathbf{x}_i'(p+1)(\mathbf{a}(h,p+1) - \breve{\mathbf{a}}_i(h,k))\}^2 \\
&= \sum_{i=m_h}^{n-h} \eta_{i,h}^2 + (1+\mathrm{o}(1)) \sum_{i=m_h}^{n-h} \{\mathbf{x}_i'(p+1)(\mathbf{a}(h,p+1) - \breve{\mathbf{a}}_i(h,k))\}^2 \\
&\quad + \mathrm{O}(1) \qquad \text{a.s.,}
\end{aligned} \tag{B.44}$$

*where the* $\breve{\mathbf{a}}_i(h,k)$ *in* (B.44) *is viewed as a* $(p+1)$*-dimensional vector with undefined entries set to 0. Direct calculations yield*

$$\begin{aligned}
\sum_{i=m_h}^{n-h} &\{\mathbf{x}_i'(p+1)(\mathbf{a}(h,p+1) - \breve{\mathbf{a}}_i(h,k))\}^2 \\
&= (\mathbf{a}(h,p+1) - \tilde{\mathbf{a}}(h,k))' V_{n-h}(k)(\mathbf{a}(h,p+1) - \tilde{\mathbf{a}}(h,k)) \\
&\quad - 2 \sum_{i=m_h}^{n-h} \mathbf{x}_i'(p+1)(\mathbf{a}(h,p+1) - \tilde{\mathbf{a}}(h,k))\mathbf{x}_i'(p+1)(\breve{\mathbf{a}}_i(h,k) - \tilde{\mathbf{a}}(h,k)) \\
&\quad + \sum_{i=m_h}^{n-h} [\mathbf{x}_i'(k)(\breve{\mathbf{a}}_i(h,k) - \tilde{\mathbf{a}}(h,k))]^2,
\end{aligned} \tag{B.45}$$

*where the* $\tilde{\mathbf{a}}(h,k)$ *in the first two terms on the right-hand side of* (B.45) *is viewed as a* $(p+1)$*-dimensional vector with undefined entries set to 0. By* (3.2) *of Lai and Wei* [15],

$$\liminf_{n \to \infty} n^{-1} V_{n-h}(k) > 0 \qquad \text{a.s.} \tag{B.46}$$

Consequently, (3.9) follows from (B.44)–(B.46), (B.42), the Cauchy–Schwarz inequality and the fact that $\mathbf{a}(h, p+1) - \bar{\mathbf{a}}(h, k) \neq \mathbf{0}$. $\qquad\square$

**Proof of Theorem 3.4.** Since $f_{2,1}(k-1) = (k-1)\sigma^2$, Theorems 3.2 and 3.3 imply

$$P(\hat{k}_{D,n}^{(1)} = p_1, \text{eventually}) = 1. \tag{B.47}$$

Applying (B.47) and Theorems 3.1–3.3, Theorem 3.4 follows. $\qquad\square$

# Appendix C

In this Appendix, we sketch the proof of Theorem 4.1. Applying an argument used in the proof of Theorem 3.5 in Wei [22], it can be shown that for $k < p_h$,

$$\liminf_{n \to \infty} \hat{\sigma}_{D,n}^2(h, k) - \hat{\sigma}_{D,n}^2(h, p_h) > 0 \qquad \text{a.s.} \tag{C.1}$$

Armed with the probability results obtained in Appendix B, one obtains for $k_1 \geq p_1$ and $k_2 \geq p_h$,

$$|\hat{\sigma}_{P,n}^2(h, k_1) - \hat{\sigma}_{D,n}^2(h, k_2)| = \mathrm{o}(\log n / n) \qquad \text{a.s.,} \tag{C.2}$$

$$|\hat{\sigma}_{P,n}^2(h, k_1) - \hat{\sigma}_{P,n}^2(h, p_1)| = \mathrm{o}(\log n / n) \qquad \text{a.s.,} \tag{C.3}$$

$$|\hat{\sigma}_{D,n}^2(h, k_2) - \hat{\sigma}_{D,n}^2(h, p_h)| = \mathrm{o}(\log n / n) \qquad \text{a.s.} \tag{C.4}$$

In addition, it can be shown that for $k_1 \geq p_1$ and $k_2 \geq p_h$,

$$\mathrm{tr}\left\{ \left( \sum_{j=k}^{n-h} \mathbf{x}_j(k)\mathbf{x}_j'(k) \right) \ddot{L}_{h,n}(k) \left( \sum_{j=k}^{n-h} \mathbf{x}_j(k)\mathbf{x}_j'(k) \right)^{-1} \ddot{L}_{h,n}'(k) \right\} \tilde{\sigma}_n^2 C_n$$

$$= \left\{ \sigma^2 \left( \sum_{j=0}^{h-1} b_j \right)^2 + f_{1,h}(k_1 - 1) \right\} C_n + \mathrm{o}(C_n) \qquad \text{a.s.} \tag{C.5}$$

and

$$\mathrm{tr}\left\{ \left( \sum_{j=k}^{n-h} \mathbf{x}_j(k)\mathbf{x}_j'(k) \right)^{-1} \left( \sum_{j=k}^{n-2h+1} \mathbf{z}_j(k)\mathbf{z}_j'(k) \right) \right\} \tilde{\sigma}_n^2 C_n$$

$$= \left\{ \sigma^2 \left( \sum_{j=0}^{h-1} b_j \right)^2 + f_{2,h}(k_2 - 1) \right\} C_n + \mathrm{o}(C_n) \qquad \text{a.s.} \tag{C.6}$$

Consequently, the asymptotic efficiency of $(\hat{O}_n, \hat{M}_n)$ follows from (C.1)–(C.6).